\documentclass[12pt]{article}
\usepackage{amsmath,amssymb,amsthm, amscd}

  \makeatletter
  \@addtoreset{equation}{section}
  \makeatother

\renewcommand{\tilde}[1]{\widetilde{#1}}
\newtheorem{dfn}{Definition}[section]
\newtheorem{thm}[dfn]{Theorem}
\newtheorem{lem}[dfn]{Lemma}
\newtheorem{prop}[dfn]{Proposition}
\newtheorem{cor}[dfn]{Corollary}

\begin{document}
\begin{center} 
{\bf {\LARGE On smoothness of minimal models of quotient singularities by finite subgroups of $SL_n(\mathbb{C})$}}
\end{center} 
\vspace{0.4cm}

\begin{center}
{\large Ryo Yamagishi}
\end{center} 
\vspace{0.4cm}

\begin{abstract}
We prove that a quotient singularity $\mathbb{C}^n/G$ by a finite subgroup $G\subset SL_n(\mathbb{C})$ has a crepant resolution only if $G$ is generated by junior elements. This is a generalization of the result of Verbitsky \cite{V}. We also give a procedure to compute the Cox ring of a minimal model of a given $\mathbb{C}^n/G$ explicitly from information of $G$. As an application, we investigate the smoothness of minimal models of some quotient singularities. Together with work of Bellamy and Schedler, this completes the classification of symplectically imprimitive quotient singularities which admit projective symplectic resolutions.
\end{abstract}

\section{Introduction}\label{1}

Crepant resolutions of singularities play key roles in various branches of algebraic geometry, and have been studied intensively. When one treats a quotient singularity by a finite group, crepant resolutions have particularly important meanings in the context of the McKay correspondence, which relates the geometry of a crepant resolution to the representation theory of the group. The aim of this article is to tackle the existence problem of crepant resolutions of quotient singularities.

Let $V$ be a finite dimensional $\mathbb{C}$-vector space and let $G\subset SL(V)$ be a finite subgroup. How can we determine the existence of a crepant resolution of the given quotient singularity $V/G$ ? When $\mathrm{dim}\;V=2$, $V/G$ is a well-known Kleinian singularity and therefore it has a unique crepant resolution. The existence of crepant resolutions for three dimensional cases is also proven (cf. \cite{R},\cite{BKR}). However, for higher dimensional cases, this is not true in general. No general criterion for existence of crepant resolutions has been known, but when $V$ is a symplectic vector space and $G$ is a subgroup of the symplectic group $Sp(V)$, there is a useful necessary condition. Verbitsky proved that $V/G$ for $G\subset Sp(V)$ admits a symplectic (or equivalently crepant) resolution only if $G$ is generated by symplectic reflections \cite{V} (see section \ref{3} for the definition).

One of our main results in this article is the following.

\begin{thm}\label{main}
 If $V/G$ for $G\subset SL(V)$ admits a crepant resolution, then $G$ is generated by junior elements.
\end{thm}

We will define an element $g\in G$ to be {\em junior} and will give a proof of the theorem in section \ref{3}. For symplectic cases, junior elements are nothing but symplectic reflections. Thus the theorem is a generalization of Verbitsky's result to nonsymplectic cases.

We also suggest a procedure to determine the (non)existence of projective crepant reslutions of $V/G$ for a given finite subgroup $G\subset SL(V)$. The idea is as follows. By a general result of birational geometry, it is known that $V/G$ always admits a minimal model \cite{BCHM}. Since crepant resolutions are nothing but smooth minimal models, it is enough to check whether each minimal model $X$ is smooth or not. To this end, we compute the Cox ring $\mathrm{Cox}(X)$, which was introduced by Hu and Keel \cite{HK}, of $X$ from $G$ (without constructing $X$ explicitly) and recover $X$ from $\mathrm{Cox}(X)$. This is done by using the similar method to one by Donten-Bury and Wi\'{s}niewski in \cite{DW}, where the authors give the Cox ring of a symplectic resolution of $V/G$ for 4-dimensional $V$ and $G$ of order 32. We generalize their method to minimal models for any finite subgroups in $SL(V)$. We give the algorithm to calculate generators of Cox rings and show several examples including Kleinian singularities in section \ref{4}. This gives a different calculation of the Cox rings of the minimal resolutions of Kleinian singularities from the ones in \cite{FGAL} and \cite{D}. Most of the calculations need a help of a computer software such as ``SINGULAR'' \cite{GPS} or ``Macaulay2'' \cite{GS}. In Appendix in section \ref{7} we give efficient ways of calculations.

In section \ref{5} we study the property of the Cox rings from the viewpoint of geometric invariant theory (GIT) and birational geometry. The spectrum $\mathfrak{X}=\mathrm{Spec}(\mathrm{Cox}(X))$ has a natural action by an algebraic torus associated to the divisor class group $\mathrm{Cl}(X)$. The crucial fact is that every minimal model can be recovered from the Cox ring as a GIT quotient of $\mathfrak{X}$ by the torus action with an appropriate stability. The sets of generic stabilities form a fan called the GIT fan on the vector space $\mathrm{Cl}(X)\otimes\mathbb{R}$. We also discuss the structure of the GIT fan. It contains some information such as the number of minimal models. We give one example in section \ref{5}. 

Let $\chi\in\mathrm{Cl}(X)\otimes\mathbb{R}$ be a stability which gives the minimal model $X$. If the (semi)stable locus $U\subset\mathfrak{X}$ associated to $\chi$ is smooth, one can check the smoothness of $X$ by looking at the torus action on $U$. However, the author does not know if this is always the case. Moreover, when the order of the group is not small, it seems almost impossible to calculate the relations of the generators of the Cox ring and  to check the smoothness of $U$ by the Jacobian criterion even if one uses a computer. Thus from these viewpoints, our method is not enough to completely answer the question raised in the second paragraph of this section. 

As an application of the description of the Cox rings, we classify all symplectically imprimitive subgroups $G\subset Sp(V)$ such that $V/G$ admits a projective crepant resolution. This was already done by Bellamy and Schedler except 6 types of groups all of which are subgroups of $Sp(4,\mathbb{C})$ \cite{BS1}. In section \ref{6} we complete the classification by studying the remaining 6 cases (Theorem \ref{classify}). It will turn out that only one group among the exceptional groups admits a crepant resolution and that this is not a new example. 

\vspace{5mm}
\noindent{\bf Acknowledgements}\\
The author is grateful to Manfred Lehn for fruitful discussions, and to Yoshinori Namikawa for helpful comments and suggesting the proof of Theorem \ref{main2}. The author also thanks Maksymilian Grab and Gwyn Bellamy for pointing out errors in the earlier versions of this article.
The author is supported by JSPS KAKENHI Grant Number JP16J04485.

\section{Quotient singularities, minimal models and Cox rings\label{2}}
Let $V$ be a complex vector space of dimension $n$ and let $G$ be a finite subgroup of $SL(V)$. Note that $G$ contains no pseudo-reflection, that is, $g\in G$ such that $\mathrm{codim}_V V^g=1$ where $V^g$ denotes the fixed subspace by $g$. It is well-known that the quotient singularity $V/G=\mathrm{Spec}\,\mathbb{C}[V]^G$ by $G$ is Gorenstein \cite{W}. Thus we can talk about the discrepancy of exceptional divisors of a birational morphism to $V/G$. 

\begin{dfn}
A {\bf minimal model} of $V/G$ is a $\mathbb{Q}$-factorial normal variety $X$ which has only terminal singularities together with a  crepant birational morphism $X\to V/G$.
\end{dfn}

Note that a nonsingular minimal model of $V/G$ is nothing but a crepant resolution of $V/G$.
Throughout this article all minimal models and crepant resolutions are assumed to be projective over $V/G$ unless otherwise stated. The following theorem is a fundamental result in birational geometry.

\begin{thm}(\cite{BCHM})
There exists a minimal model $X$ of $V/G$. 
\end{thm}

To introduce Cox rings, we should consider the divisor class group of a minimal model. As for $V/G$, its divisor class group $\mathrm{Cl}(V/G)$ is known.

\begin{prop}\label{torsion}(\cite[Ch. 3]{Ben})
The divisor class group $\mathrm{Cl}(V/G)$ of the quotient singularity $V/G$ is canonically isomorphic to the group $Ab(G)^\vee=\mathrm{Hom}(G, \mathbb{C}^*)$ of characters of $G$ where $Ab(G)$ denotes the abelianization $G/[G,G]$ of G. In particular $\mathrm{Cl}(V/G)$ is a torsion group.
\end{prop} 

Let $X(\stackrel{\pi}{\to}V/G)$ be a minimal model. Then one easily sees that $\mathrm{Cl}(X)$ is also finitely generated since every divisor of $X$ consists of exceptional divisors of $\pi$ and divisors from $V/G$. We now assume that $\mathrm{Cl}(X)$ is torsion free for simplicity. Let $D_1,\dots,D_m$ be Weil divisors whose classes form a basis of $\mathrm{Cl}(X)$. We define the Cox ring of $X$ as
$$\mathrm{Cox}(X)=\bigoplus_{(a_1,\dots,a_m)\in\mathbb{Z}^m}H^0(X,\mathcal{O}_X(a_1D_1+\cdots+a_mD_m)).$$
For a Weil divisor $D$, the vector space $H^0(X,\mathcal{O}_X(D))$ is identified with the set
$$\{f\in\mathbb{C}(X)^*|\mathrm{div}(f)+D\ge0\}\cup\{0\},$$
and the Cox ring has the natural ring structure inherited from the multiplication in $\mathbb{C}(X)$. It is known that the isomorphism class of the Cox ring is independent of the choice of $D_i$'s (cf. \cite{ADHL}). $\mathrm{Cox}(X)$ has a $\mathrm{Cl}(X)$-grading and this grading gives a torus action on $\mathrm{Cox}(X)$ in the following way. Let $T:=\mathrm{Hom(\mathrm{Cl}(X),\mathbb{C}^*)}$ be the algebraic torus. It acts on the homogeneous part $H^0(X,\mathcal{O}_X(a_1D_1+\cdots+a_mD_m))$ of $\mathrm{Cox}(X)$ for $D=\sum_{i=1}^ma_iD_i\in\mathrm{Cl}(X)$ by multiplying $t(D)$ for each $t\in T$. This action naturally induces an action on the spectrum $\mathfrak{X}=\mathrm{Spec(\mathrm{Cox}(X))}$.

Next we consider GIT quotients of $\mathfrak{X}$ by $T$. To this end, we should choose a $T$-linearization on $\mathfrak{X}$. We particularly use the trivial line bundle twisted by a character of $T$. When we take a divisor class $D\in\mathrm{Cl}(X)$, we can regard it as a character of $T$ by the evaluation map $T\to\mathbb{C}^*,t\mapsto t(D)$. We can define the GIT quotient $\mathfrak{X}/\!/_D T$ of $\mathfrak{X}$ by $T$ with respect to a character $D$ of $T$. As we will see later in section \ref{5}, the most important feature of $\mathfrak{X}$ is that every minimal model $X'$ of $V/G$ can be obtained as a GIT quotient of $\mathfrak{X}$ for some $D$.


\section{Discrete valuations on function fields\label{3}}
Let $V$ and $G$ be as in the previous section. For an element $g\in G$, we define a discrete valuation $\nu_g:\mathbb{C}(V)\to\mathbb{Z}\cup\{\infty\}$ on the rational function field $\mathbb{C}(V)$ of $V$ as follows. Let $x_1,\dots,x_n\in V^*$ be linearly independent eigenvectors of $g$. Then there are unique integers $a_i$ for $i=1,\dots,n$ such that  $0\le a_i<r$ and
$$g\cdot x_i=\mathrm{exp}\,\frac{2\pi\sqrt{-1}a_i}{r}\cdot x_i$$
where $r$ is the order of $g$. For any nonzero polynomial
$$f=\sum_{\alpha=(\alpha_1,\dots,\alpha_n)\in\mathbb{Z}_{\ge0}^n}c_\alpha x_1^{\alpha_1}\cdots x_n^{\alpha_n},\;c_\alpha\in\mathbb{C},$$
we set
$$\nu_g(f)=\underset{{\alpha=(\alpha_1,\dots,\alpha_n)\atop c_\alpha\ne0}}{\mathrm{min}}\{\sum_i \alpha_i a_i\}.$$
This extends uniquely to a discrete valuation on the whole of $\mathbb{C}(V)$.

If $\pi:X\to V/G$ is a minimal model, the function field $\mathbb{C}(X)$ of $X$ is identified with the $G$-invariant subfield $\mathbb{C}(V)^G$ of $\mathbb{C}(V)$. Therefore we get a discrete valuation on $\mathbb{C}(X)$ by restricting $\nu_g$. On the other hand, we also have another valuation on $\mathbb{C}(X)$. Let $E$ be an irreducible exceptional divisor of $\pi$. Then it gives the divisorial valuation $\nu_E$ on $\mathbb{C}(X)$ defined by $\nu_E(f)=\mathrm{ord}_E(\mathrm{div}(f))$ for $f\in\mathbb{C}(X)^*$.

Next we introduce the notion of the {\it age} of an element $g\in G$. Let $a_1,\dots,a_n$ and $r$ be as above. Then we set $\mathrm{age}(g)=\frac{1}{r}\sum_{i=1}^n a_i$. Note that $\mathrm{age}(g)$ is always an integer since $g$ is in $SL(V)$. One can easily check that age is invariant under conjugation by $GL(V)$.

\vspace{3mm}
{\bf Remark.}\;In \cite{IR}, age is not defined as a function on $G$ but on $\Gamma:=\mathrm{Hom}(\mu_R,G)$ where $R\in \mathbb{N}$ is a common multiple of the orders of all elements in $G$ and $\mu_R$ is the group of the $R$-th roots of unity. However, the definition in \cite{IR} coincides with ours via the isomorphism $\Gamma\to G,f\mapsto f(\mathrm{exp}\,\frac{2\pi\sqrt{-1}}{R})$.

\vspace{3mm}
We call an element $g\in G$ {\em junior} if $\mathrm{age}(g)=1$. The following theorem claims that the information about exceptional divisors of a minimal model can be read off from the information about $G$. Ito and Reid proved the following.

\begin{thm}\label{IR}(\cite{IR})
Let $\pi:X\to V/G$ be a proper birational morphism from a $\mathbb{Q}$-factorial normal variety $X$. Then $\pi$ is a (not necessarily projective) minimal model if and only if there is a bijection of the sets
$$\{\mbox{an irreducible exceptional divisor of }\pi\}\cong\{\mbox{a conjugacy class of junior elements in G}\}$$
such that, if $E$ is an irreducible exceptional divisor of $\pi$ which corresponds to $g\in G$ via this bijection, the equality
$$\nu_E=\frac{1}{r}\nu_g|_{\mathbb{C}(X)}$$
holds.
\end{thm}

We also give a result in relation to the existence of a smooth minimal model (i.e. a crepant resolution) of $V/G$.

\begin{dfn}
Let $V$ be a finite dimensional symplectic $\mathbb{C}$-vector space and let $G$ be a finite subgroup of $Sp(V)$. An element $g\in G$ is called a {\bf symplectic reflection} if the codimension of the fixed point set $V^g$ in $V$ is two.
\end{dfn} 

Verbitsky proved that if $V/G$ admits a (not necessarily projective) symplectic (or equivalently crepant) resolution, then $G$ is generated by symplectic reflections \cite{V}. The aim of this section is to generalize this result.

\begin{thm}\label{main2}
Let $G$ be a finite subgroup of $SL(V)$ for a finite dimensional $\mathbb{C}$-vector space $V$ and let $\pi:X\to V/G$ be a minimal model. Then the algebraic fundamental group $\pi_1^{\mathrm{alg}}(X_{\mathrm{reg}})$ of the regular part of $X$ is trivial if and only if $G$ is generated by junior elements.
\end{thm}

{\em Proof.} Let $H\subset G$ be the normal subgroup generated by all junior elements in $G$ and let $X'$ (resp. $X''$) be the main component (which dominates $V/G$) of the normalization of the fibred product $X\times_{V/G}V/H$ (resp. $X\times_{V/G}V$). Then we have the commutative diagram
$$\begin{CD}X'' @>\tilde{p}'>> X' @>\tilde{p}>> X\\
@V\pi'' VV @V\pi' VV @VV\pi V\\
V@>p'>> V/H@>p>> V/G.\\
\end{CD}$$

We need the following lemma.



\vspace{3mm}
\begin{lem}\label{unramify}
Let $E$ be a $\pi$-exceptional divisor. Then the finite surjective morphism $\tilde{p}$ is unramified at any generic point of $\tilde{p}^{-1}(E)$.
\end{lem}

{\em Proof.} Let $F_1,\dots, F_k$ be the irreducible components of $\tilde{p}^{-1}(E)$ and let $F'_i$ be an irreducible component of $\tilde{p}'^{-1}(F_i)$. By construction, one has the equality of valuations on $\mathbb{C}(X)$
$$\nu_{F'_i}|_{\mathbb{C}(X)}=r_1\nu_{F_i}|_{\mathbb{C}(X)}=r_2\nu_E$$
for any $i$ where $r_1$ and $r_2$ are the ramification indices along $F'_i$ of $\tilde{p}'$ and $\tilde{p}\circ\tilde{p}'$ respectively. Let $g\in G$ be an element in the conjugacy class corresponding to $E$ via the bijection in Theorem \ref{IR}. By \cite[2.6 and 2.8]{IR}, one has $r_1=r_2=\sharp\langle g\rangle$. Therefore the claim holds.
\qed

\vspace{3mm}
{\bf Remark.}\;Since $\nu_{F_i}|_{\mathbb{C}(X)}=\nu_E$, one can check that $\pi':X'\to V/H$ is a minimal model by Theorem \ref{IR}. Let $C\subset G$ be the $G$-conjugacy class containing $g$. Then the decomposition of $\tilde{p}^{-1}(E)$ into the irreducible components corresponds to the division of $C$ into the $H$-conjugacy classes.

\vspace{3mm}
Now we return to the proof of the theorem. First we assume that $G\ne H$. By this assumption and Lemma \ref{unramify}, the map $\tilde{p}$ is \'{e}tale in codimension one of $\deg(\tilde{p})>1$. By the purity of branch locus, this implies that $\tilde{p}$ is \'{e}tale over $X_{\mathrm{reg}}$. Therefore $\pi_1^{\mathrm{alg}}(X_{\mathrm{reg}})$ is nontrivial. 

Conversely, we assume that $\pi_1^{\mathrm{alg}}(X_{\mathrm{reg}})\ne 1$. Then there is a nontrivial finite \'{e}tale covering $Y_0\to X_{\mathrm{reg}}$. By taking the normalization of $X$ in $\mathbb{C}(Y_0)$, one can extend the covering map to a finite surjective morphism $q:Y\to X$. Let $Y\to Z\to V/G$ be the Stein factorization of $\pi\circ q$. As $Z\to V/G$ is finite \'{e}tale over $(V/G)_{\mathrm{reg}}$, we can write $Z=V/K$ for a suitable subgroup $K$ of $G$. Since $q$ is \'{e}tale in codimension one, the birational morphism $Y\to V/K$ is also a minimal model. By Theorem \ref{IR}, we see that $K$ contains all junior elements. (Any junior element $g\in G$ defines an exceptional divisor of $Y\to V/K$, and thus $g$ is $G$-conjugate to an element in $K$.) Therefore $H\subset K\ne G$.
\qed

\vspace{3mm}
Since symplectic reflections are nothing but junior elements (see e.g. \cite{Ka1} Lemma 1.1), the following is a generalization of Verbitsky's result.

\begin{cor}(=Theorem \ref{main})
If $V/G$ admits a (not necessarily projective) crepant resolution, then $G$ is generated by junior elements.
\end{cor}

{\em Proof.} If $X\to V/G$ is a crepant resolution, the fundamental group $\pi_1(X)$ is trivial (see \cite[Thm. 7.8]{Ko} or \cite[ Thm. 4.1]{V}). As $X=X_{\mathrm{reg}}$, we conclude by the theorem that $G$ is generated by junior elements. Note that we did not use the projectivity of $X\to V/G$.
\qed

\section{Embedding of the Cox ring and description of the generators\label{4}}
The goal of this section is to give  an explicit procedure for calculating the Cox ring of a minimal model of a given $V/G$. This is done by considering the Cox ring as a subring of some bigger and simpler ring. This construction is almost due to Donten-Bury and Wi\'{s}niewski. In \cite{DW} the authors calculated the Cox ring for a group of order 32 acting on a 4-dimensional symplectic vector space (see Example 2 below in this section). 

As in section \ref{2}, we can also define the Cox ring of $V/G$. It is defined as
$$\mathrm{Cox}(V/G)=\bigoplus_{D\in \mathrm{Cl}(V/G)}H^0(V/G, \mathcal{O}_{V/G}(D))$$
as a $\mathbb{C}$-vector space where $\mathcal{O}_{V/G}(D)$ is the rank-1 reflexive sheaf associated to a Weil divisor class $D$. Note that $H^0(V/G, \mathcal{O}_{V/G}(D))$ is identified with $\{f\in\mathbb{C}(V/G)^*|\mathrm{div}(f)+D'\ge0\}\cup\{0\}$ where $D'$ is any Weil divisor on $V/G$ which represents $D$. Then $\mathrm{Cox}(V/G)$ has a $\mathrm{Cl}(V/G)$-graded ring structure which is defined similarly to the case in section \ref{2}. However, this construction is not exactly the same because of torsions in $\mathrm{Cl}(V/G)$ (cf. Proposition \ref{torsion}). See \cite{ADHL} for details.

The degree zero part of $\mathrm{Cox}(V/G)$ is $\mathbb{C}[V]^G$, and thus $\mathrm{Cox}(V/G)$ is a $\mathbb{C}[V]^G$-algebra. We have the following result.

\begin{prop}(\cite[Theorem 3.1]{AG}\label{AG})
There is an isomorphism as $\mathbb{C}[V]^G$-algebras between $\mathrm{Cox}(V/G)$ and $\mathbb{C}[V]^{[G,G]}$ which preserves the natural gradings by $Ab(G)^\vee$.
\end{prop} 

Let $g_1,\dots,g_m$ be the complete system of representatives of the conjugacy classes of junior elements in $G$ and set $\nu_i:=\nu_{g_i}$. Then for each $i$ there is a unique irreducible exceptional divisor $E_i$ of $\pi:X\to V/G$ such that $\nu_{E_i}=\frac{1}{r}\nu_i|_{\mathbb{C}(X)}$ by Theorem \ref{IR}. Let $\mathrm{Cl}(X)^{\mathrm{free}}$ denote the free abelian group $\mathrm{Cl}(X)/\mathrm{Cl}(X)^{\mathrm{tor}}$ where $\mathrm{Cl}(X)^{\mathrm{tor}}$ is the torsion part of $\mathrm{Cl}(X)$. Then the rank of $\mathrm{Cl}(X)^{\mathrm{free}}$ is $m$. This follows from the short exact sequence
\begin{equation}\label{exact}
0\to\bigoplus_{i=1}^m \mathbb{Z}E_i\to \mathrm{Cl}(X)\stackrel{\pi_*}{\to}\mathrm{Cl}(V/G)\to 0
\end{equation}
noticing that $\mathrm{Cl}(V/G)$ is a torsion group (Proposition \ref{torsion}).

Let $\mathbb{C}[\mathrm{Cl}(X)^{\mathrm{free}}]=\bigoplus_{\bar{D}\in \mathrm{Cl}(X)^{\mathrm{free}}}\mathbb{C}t^{\bar{D}}$ be the group algebra where $t^{\bar{D}}$'s denote the basis. Now we construct an embedding of $\mathrm{Cox}(X)$ into $\mathbb{C}[V]^{[G,G]}\otimes_{\mathbb{C}}\mathbb{C}[\mathrm{Cl}(X)^{\mathrm{free}}]$ as follows. For any Weil divisor $D$ on $X$ and any homogeneous element $\tilde{f}\in H^0(X,\mathcal{O}_X(D))\subset\mathrm{Cox}(X)$, we can regard $\tilde{f}$ as an element of $H^0(V/G, \mathcal{O}_{V/G}(\pi_*D))$ via the identification $\mathbb{C}(X)=\mathbb{C}(V/G)$ between the function fields. Let $f\in \mathbb{C}[V]^{[G,G]}$ be the corresponding element to $\tilde{f}$ via the isomorphism appeared in Proposition 4.1. Then we obtain a ring homomorphism $\Theta:\mathrm{Cox}(X)\to \mathbb{C}[V]^{[G,G]}\otimes_{\mathbb{C}}\mathbb{C}[\mathrm{Cl}(X)^{\mathrm{free}}]$ setting $\Theta(\tilde{f})=f\otimes t^{\bar{D}}$ where $\bar{D}$ is the class of $D$ in $\mathrm{Cl}(X)^{\mathrm{free}}$. The following is a generalization of \cite[Prop. 3.8]{DW}.

\begin{lem}
$\Theta:\mathrm{Cox}(X)\to \mathbb{C}[V]^{[G,G]}\otimes_{\mathbb{C}}\mathbb{C}[\mathrm{Cl}(X)^{\mathrm{free}}]$ is injective.
\end{lem}

{\em Proof.} Let $\tilde{f}$ be any element in the kernel of $\Theta$. As $\Theta$ is compatible with the quotient map $\mathrm{Cl}(X)\to\mathrm{Cl}(X)^{\mathrm{free}}$, we may assume that all the divisor classes $D_i$'s to which the homogeneous components of $\tilde{f}$ belong are in the same class of $\mathrm{Cl}(X)^{\mathrm{free}}$. On the other hand, the natural map $\mathrm{Cox}(X)\to\mathrm{Cox}(V/G)$ which is obtained by the composition of $\Theta$ followed by the evaluation $t=1$ is also compatible with the surjection $\mathrm{Cl}(X)\to\mathrm{Cl}(V/G)$. Therefore  we may also assume that $D_i$'s are in the same class of $\mathrm{Cl}(X)/\bigoplus_{i=1}^r \mathbb{Z}E_i$ by (\ref{exact}). However, since the subgroup $\bigoplus_{i=1}^r \mathbb{Z}E_i$ is torsion-free, the element $\tilde{f}\in\mathrm{Cox}(X)$ must be homogeneous. In this case the claim is clear by definition.
\qed

\vspace{5mm}

Therefore $\mathrm{Cox}(X)$ can be realized as a subring of $\mathbb{C}[V]^{[G,G]}\otimes_{\mathbb{C}}\mathbb{C}[\mathrm{Cl}(X)^{\mathrm{free}}]$. Our next task is to know which elements in $\mathbb{C}[V]^{[G,G]}\otimes_{\mathbb{C}}\mathbb{C}[\mathrm{Cl}(X)^{\mathrm{free}}]$ are in the image of $\Theta$.

Let $p:V/[G,G]\to V/G$ be the quotient map. For an element $f$ of $\mathbb{C}[V]^{[G,G]}$ which is homogeneous with respect to $Ab(G)^\vee$-grading, consider the Weil divisor $D_f=p_*(\mathrm{div}_{V/[G,G]}(f))$ on $V/G$. Let $\bar{D}_f$ be the class in $\mathrm{Cl}(X)^{\mathrm{free}}$ of the strict transform of $D_f$ by $\pi^{-1}:V/G\dashrightarrow X$.

\begin{lem}\label{lem:strict transform}
Let $f$ and $\bar{D}_f$ be as above. Then the equality
$$\bar{D}_f=-\sum_{i=1}^m\frac{1}{r_i} \nu_i(f)\bar{E}_i$$
in $\mathrm{Cl}(X)^{\mathrm{free}}$ holds where $\bar{E}_i$ is the class of $E_i$ in $\mathrm{Cl}(X)^{\mathrm{free}}$.
Moreover, $f\otimes t^{\bar{D}_f}$ is in $\mathrm{Im}\,\Theta$.
\end{lem}

{\em Proof.} Let $\tilde{D}_f$ be the strict transform of $D_f$ by $\pi^{-1}:V/G\dashrightarrow X$ and let $\tilde{f}$ be the element of $\mathrm{Cox}(V/G)$ which corresponds to $f$ via the isomorphism in Proposition \ref{AG}. As $f$ is homogeneous, some power $f^r$ ($r\in\mathbb{N}$) is in $\mathbb{C}[V]^G$ and  equals $\tilde{f}^r\in\mathrm{Cox}(V/G)_0$. Since the pullback of $rD_f$ by $\pi$ can be written as
$$\pi^*(\mathrm{div}_{V/G}(\tilde{f}^r))=r\tilde{D}_f+\nu_{E_1}(\tilde{f}^r)E_1+\cdots+\nu_{E_m}(\tilde{f}^r)E_m,$$
we have $r\tilde{D}_f=-r\sum_{i=1}^m\frac{1}{r_i} \nu_i(f)E_i$ in $\mathrm{Cl}(X)$ by Theorem \ref{IR}. By dividing the both sides by $r$, we obtain the desired equation.

For the second claim, one can easily check by definition that $\Theta(\tilde{f})=f\otimes t^{\bar{D}_f}$.
\qed

\vspace{5mm}

By this lemma, we can describe generators of $\mathrm{Im}\,\Theta$ from those of $\mathbb{C}[V]^{[G,G]}$. Let $S=\{\phi_1,\dots,\phi_k\}$ be a generating system of $\mathbb{C}[V]^{[G,G]}$ such that each $\phi_j$ is homogeneous with respect to the $Ab(G)$-action. We consider the following condition ($*$):\\
``For every nonzero $f\in\mathbb{C}[V]^{[G,G]}$, there are monomials $f_1,\dots,f_l$ of $\phi_1,\dots,\phi_k$ such that\\
(i) $f=f_1+\cdots+f_l$, and\\
(ii) $\nu_i(f)\le\nu_i(f_j)$ for every $i$ and every $j$.''

\begin{prop}\label{prop:generator}
If the homogeneous generators $\phi_1,\dots,\phi_k$ of $\mathbb{C}[V]^{[G,G]}$ satisfy ($*$), then the subset
$$\{\phi_j\otimes t^{\bar{D}_{\phi_j}}\}_{j=1,\dots,k}\cup\{t^{\bar{E}_1},\dots,t^{\bar{E}_m}\}$$
of $\mathbb{C}[V]^{[G,G]}\otimes_{\mathbb{C}}\mathbb{C}[\mathrm{Cl}(X)^{\mathrm{free}}]$ is a generating system of $\mathrm{Im}\,\Theta$.
\end{prop}

{\em Proof.} First note that $t^{\bar{E}_1},\dots,t^{\bar{E}_m}$ are in $\mathrm{Im}\,\Theta$ since $E_i$'s are effective divisors. Take any homogeneous element $f\otimes t^{\bar{D}}$ in $\mathrm{Im}\,\Theta$ with a Weil divisor $D$ on $X$. By the condition ($*$) we can write $f=f_1+\cdots+f_l$ satisfying the conditions and hence can also write
$$f\otimes t^{-\sum_{i=1}^m\frac{1}{r_i} \nu_i(f)\bar{E}_i}=\sum_{j=1}^l f_j\otimes t^{-\sum_{i=1}^m(\frac{1}{r_i} \nu_i(f_j)\bar{E}_i-\sum_i a_{i,j}\bar{E}_i)}$$
for some $a_{i,j}\in\mathbb{Q}_{\ge 0}$. Since the images of the RHS and $\sum_{j=1}^l f_j\otimes t^{-\sum_{i=1}^m\frac{1}{r_i} \nu_i(f_j)\bar{E}_i}$ by the natural map $\mathrm{Im}\,\Theta\to\mathbb{C}[V]^{[G,G]}$ are the same, the sum $\sum_i a_{i,j}\bar{E}_i$ must be in $\bigoplus_{i=1}^m\mathbb{Z}_{\ge0}\bar{E}_i$.

Take $\tilde{f}$ so that $\Theta(\tilde{f})=f\otimes t^{\bar{D}}$. Then the inequality $\sum_{i=1}^m\frac{1}{r_i} \nu_i(f)\bar{E}_i+D\ge0$ must be satisfied since $\tilde{f}$ is in $H^0(X,\mathcal{O}(D))$. Thus we can write $f\otimes t^{\bar{D}}=f\otimes t^{-\sum_{i=1}^m\frac{1}{r_i} \nu_i(f)\bar{E}_i}\cdot  t^{\sum_i b_{i,j}\bar{E}_i}$ for some $b_{i,j}\in\mathbb{Q}_{\ge 0}$. The same argument as above shows that $b_{i,j}$ is in $\mathbb{Z}_{\ge0}$. Since each $f_j$ is a monomial of $\phi_1,\dots,\phi_k$, each $f_j\otimes t^{\bar{D}_{f_j}}=f_j\otimes t^{-\sum_{i=1}^m\frac{1}{r_i} \nu_i(f_j)\bar{E}_i}$ is also a monomial of 
$\phi_j\otimes t^{\bar{D}_{\phi_j}}$'s. Therefore we obtain a desired expression of $f\otimes t^{\bar{D}}$.
\qed

\vspace{5mm}

Therefore we can construct the Cox ring explicitly if we find generators of $\mathbb{C}[V]^{[G,G]}$ satisfying ($*$). Now we give an algorithm for finding such generators from any generators of $\mathbb{C}[V]^{[G,G]}$.

Let $\phi_1,\dots,\phi_k$ be generators of $\mathbb{C}[V]^{[G,G]}$. We may assume that they are homogeneous with respect to the $Ab(G)$-action. For each $i\in\{1,\dots,m\}$, let $x_{i,1},\dots,x_{i,n}\in V^*$ be linearly independent eigenvectors for the $\langle g_i \rangle$-action. When we write an element $\phi\in\mathbb{C}[V]^{[G,G]}$ as the sum of monomials of $x_{i,1},\dots,x_{i,n}$, let $\mbox{min}_i(\phi)$ be the sum of the monomials whose values of $\nu_i$ are minimal among these monomials.

Consider the ring homomorphism $\alpha:\mathbb{C}[X_1,\dots,X_k]\to\mathbb{C}[V]^{[G,G]},\;X_j\mapsto\phi_j$ and let $I$ be the kernel of $\alpha$. We give a grading on $\mathbb{C}[X_1,\dots,X_k]$ by setting $\deg_i(X_j)=\nu_i(\phi_j)$. For an inhomogeneous polynomial $h$, we promise to let $\deg_i(h)$ denote the {\bf minimal} degree of $h$. We can define the minimal part $\mbox{min}_i(h)$ for $h\in\mathbb{C}[X_1,\dots,X_k]$, and let $\mbox{min}_i(I)\subset\mathbb{C}[X_1,\dots,X_k]$ be the ideal generated by the set $\{\mbox{min}_i(h)|h\in I\}$. One sees that, for each $h\in\mbox{min}_i(I)$, there is $\tilde{h}\in \mathbb{C}[X_1,\dots,X_k]$ such that $h-\tilde{h}\in I$ and $\deg_i(h)<\deg_i(\tilde{h})$.

On the other hand, consider the ring homomorphism $\beta_i:\mathbb{C}[X_1,\dots,X_k]\to\mathbb{C}[V],\;X_j\mapsto \mbox{min}_i(\phi_j)$ and let $J$ be the kernel of $\beta_i$. Let $\mbox{min}_i(J)$ denote the ideal generated by homogeneous elements in $J$ with respect to the $Ab(G)$-action. (Note that the $Ab(G)$-action on $\mathbb{C}[V]^{[G,G]}$ naturally lift to $\mathbb{C}[X_1,\dots,X_k]$.) Thus we obtain two ideals of the polynomial ring associated to $\phi_1,\dots,\phi_k$: the ideal $\mbox{min}_i(I)$ of the minimal terms of the relations, and the ideal $\mbox{min}_i(J)$ of the $Ab(G)^\vee$-homogeneous relations of the minimal terms. The following lemma is straightforward.

\begin{lem}
The inclusion $\mathrm{min}_i(I)\subset\mathrm{min}_i(J)$ holds.
\end{lem}



For a subset $A\subset\{1,\dots,m\}$, let $R_A$ be the polynomial ring $\mathbb{C}[X_1,\dots,X_k,\{t_i\}_{i\in A}]$. The grading $\deg_i$ ($i\in A$) on $\mathbb{C}[X_1,\dots,X_k]$ naturally extends to one on $R_A$ by setting $\deg_i(t_j)=-\delta_{i,j}$ where $\delta_{i,j}$ is Kronecker delta. For nonzero $h\in \mathbb{C}[X_1,\dots,X_k]$, let $h_A$ be the element in $R_A$ obtained by homogenizing $h$ by $t_i$'s with respect to $\deg_i$'s respectively. That is, if $h=\sum_l h_l$ where $h_l$'s are monomials, then $h_A=\sum_l (h_l \prod_{i\in A} t_i^{\deg_i(h_l)-\deg_i(h)})$. Consider the ideal $I=\mathrm{Ker}\,\alpha$ as above and let $I_A\subset R_A$ be the homogeneous ideal generated by the set $\{h_A|h\in I\}$.

Now we define a collection $\{S_p\}_{p=0,1,\dots}$ of subsets each of which consists of finitely many elements in $\mathbb{C}[V]^{[G,G]}$ inductively by taking the following steps.

\vspace{3mm}
\noindent{\bf Step (0)}\\
Set $S=S_0:=\{\phi_1,\dots,\phi_k\}$ and go to {\bf Step ($\boldsymbol{0,1}$)}.

\vspace{3mm}
\noindent{\bf Step ($\boldsymbol{0,i}$)} ($i=1,\dots,m$)\\
Compute $\mathrm{min}_i(I)$ and $\mathrm{min}_i(J)$ for $S$. Write $\mathrm{min}_i(J)=\mathrm{min}_i(I)+(h_1,\dots,h_l)$ with $\deg_i$-homogeneous $h_j\not\in\mathrm{min}_i(I)$. (Note that $\mathrm{min}_i(I)$ and $\mathrm{min}_i(J)$ are $\deg_i$-homogeneous ideals.)\\
$\bullet$ If $\mathrm{min}_i(I)=\mathrm{min}_i(J)$ and $i=m=1$, set $S_{p+1}=S_{p+2}=\cdots:=S_p$.\\
$\bullet$ If $\mathrm{min}_i(I)=\mathrm{min}_i(J)$ and $i=1<m$, replace $S$ by $S_{p+1}:=S_p$ and go to {\bf Step ($\boldsymbol{0,2}$)}.\\
$\bullet$ If $\mathrm{min}_i(I)=\mathrm{min}_i(J)$ and $i>1$, replace $S$ by $S_{p+1}:=S_p$ and go to {\bf Step ($\boldsymbol{1,i}$)}.\\
$\bullet$ If $\mathrm{min}_i(I)\subsetneq\mathrm{min}_i(J)$, replace $S$ by $S_{p+1}:=S_p\cup\{\alpha(h_1),\dots,\alpha(h_l)\}$ and go to {\bf Step ($\boldsymbol{0,i}$)} again.

\vspace{3mm}
\noindent{\bf Step ($\boldsymbol{i',i}$)}($1\le i'<i\le m$)\\
Compute the two ideals $\tilde{I}_{i',i}:=I_{\{1,\dots,i',i\}}\cap(t_{i'},t_i)$ and $\tilde{I}'_{i',i}:=(I_{\{1,\dots,i',i\}}\cap(t_{i'}))+(I_{\{1,\dots,i',i\}}\cap(t_i))$ for $S$. Write
$$\tilde{I}_{i',i}=\tilde{I}'_{i',i}+(h_1,\dots,h_l)$$
where $h_j$'s are elements in $R_{\{1,\dots,i',i\}}\backslash\tilde{I}'_{i',i}$ which are homogeneous with respect to $\deg_j$ for each $j\in\{1,\dots,i',i\}$. (Note that $\tilde{I}_{i',i}$ and $\tilde{I}'_{i',i}$ are $\deg_j$-homogeneous ideals for each $j\in\{1,\dots,i',i\}$.)\\
$\bullet$ If $\tilde{I}_{i',i}=\tilde{I}'_{i',i}$, $i=m$ and $i'=m-1$, then set $S_{p+1}=S_{p+2}=\cdots:=S_p$.\\
$\bullet$ If $\tilde{I}_{i',i}=\tilde{I}'_{i',i}$, $i<m$ and $i'=i-1$, then replace $S$ by $S_{p+1}:=S_p$ and go to {\bf Step ($\boldsymbol{0,i+1}$)}.\\
$\bullet$ If $\tilde{I}_{i',i}=\tilde{I}'_{i',i}$, $i<m$ and $i'<i-1$, then replace $S$ by $S_{p+1}:=S_p$ and go to {\bf Step ($\boldsymbol{i'+1,i}$)}.\\
$\bullet$ If $\tilde{I}_{i',i}\supsetneq\tilde{I}'_{i',i}$, then replace $S$ by $S_{p+1}:=S_p\cup\{\alpha(\mathrm{min}_i(h_1|_{t=1})),\dots,\alpha(\mathrm{min}_i(h_l|_{t=1}))\}$ and go to {\bf Step ($\boldsymbol{i',i}$)} again.

\vspace{3mm}

We should perform the above algorithm in the following order of the steps:
\begin{center}
$(0)\to(0,1)\to(0,2)\to(1,2)
\to(0,3)\to(1,3)\to(2,3)\to(0,4)\to(1,4)\to\cdots$
\end{center}

Note that each $S_p$ is not unique since it involves several choices. A concrete procedure for performing the algorithm above (with a computer) will be given in Appendix in section \ref{7}.

We will show that the algorithm gives generators of the Cox ring  (cf. Corollary \ref{cor:finite}). To this end, we define conditions on $S$ as follows.

\begin{dfn}
Fix $S=S_p$.\\
$\bullet$\;Let $A\in\{1,\dots,m\}$ be any subset. We say that $S$ satisfies ($*A$) if for every nonzero $f\in\mathbb{C}[V]^{[G,G]}$, there is $h\in\mathbb{C}[X_1,\dots,X_k]$ such that $f=\alpha(h)$ and $\nu_i(f)=\deg_i(h)$ for every $i\in A$.\\
$\bullet$\;Let $A\subset\{1,\dots,m\}\backslash\{i\}$ be any subset and let $h\in\mathbb{C}[X_1,\dots,X_k]$ be any element. We say that $h$ satisfies ($*A,i,p$) if there is $\tilde{h}\in\mathbb{C}[X_1,\dots,X_k]$ such that $h-\tilde{h}\in I$, $\deg_i(h)<\deg_i(\tilde{h})$, and $\deg_j(h)\le\deg_j(\tilde{h})$ for every $j\in A$.
\end{dfn}

Note that the inequality $\nu_i(f)\ge\deg_i(h)$ always holds. However, the equality does not hold in general. 

It is known that the Cox ring of a minimal model $X$ is finitely generated \cite{BCHM}. Therefore we can take finitely many homogeneous generators $\tilde{f}_1,\dots,\tilde{f}_s$ of $\mathrm{Im}\,\Theta$. Set $f_j:=\tilde{f}_j|_{t=1}\in\mathbb{C}[V]^{[G,G]}$.

\begin{lem}\label{lem:reduction}
Fix $S=S_p$ and let $A\subset\{1,\dots,m\}\backslash\{i\}$ be any subset. Assume that $S_p$ satisfies ($*A$). Then $S_p$ satisfies ($*A\cup\{i\}$) if and only if $h$ satisfies ($*A,i,p$) for any $\deg_i$-homogeneous element $h\in\mathrm{min}_i(J)$.
\end{lem}

{\em Proof.} We first assume that $S_p$ satisfies ($*A\cup\{i\}$). Take any $\deg_i$-homogeneous element $h$ from $\mathrm{min}_i(J)$. Then one sees that $\nu_i(\alpha(h))>\deg_i(h)$. On the other hand, there is $\tilde{h}\in\mathbb{C}[X_1,\dots,X_k]$ such that $\alpha(\tilde{h})=\alpha(h)$ and $\nu_j(\alpha(\tilde{h}))=\deg_j(\tilde{h})$ for every $j\in A\cup\{i\}$ since $S_p$ satisfies ($*A\cup\{i\}$). Therefore $h$ satisfies ($*A,i,p$).

To prove the converse, we show that $S_p$ satisfies ($*A\cup\{i\}$) for each $f\in\{f_1,\dots,f_s\}$ assuming that $h$ satisfies ($*A,i,p$) for any $\deg_i$-homogeneous $h\in\mathrm{min}_i(J)$. By assumption there is $h'\in\mathbb{C}[X_1,\dots,X_k]$ such that $f=\alpha(h')$ and $\nu_j(f)=\deg_j(h')$ for every $j\in A$. If $\nu_i(f)=\deg_i(h')$, we are done. So we assume otherwise. Then the $\deg_i$-minimal part $h=\mathrm{min}_i(h')$ is in $\mathrm{min}_i(J)$. Since $h$ satisfies ($*A,i,p$), there is $\tilde{h}\in\mathbb{C}[X_1,\dots,X_k]$ such that $h-\tilde{h}\in I$, $\deg_i(h)<\deg_i(\tilde{h})$, and $\deg_j(h)\le\deg_j(\tilde{h})$ for every $j\in A$. Therefore replacing $h$ in $h'$ by $\tilde{h}$ increases the value of $\deg_i$ without decreasing the values of $\deg_j$ for $j\in A$. Repeating this process gives $h'$ such that $\nu_i(f)=\deg_i(h')$ by induction on $\nu_i(f)-\deg_i(h')$.

By the argument in the proof of Proposition \ref{prop:generator}, we see that each $\tilde{f}_j\in\mathrm{Im}\,\Theta$ is a product of  $f_j\otimes t^{\bar{D}_{f_j}}$ and $t^{\bar{E}_i}$'s. Therefore, for any nonzero $f\in\mathbb{C}[V]^{[G,G]}$, the element $f\otimes t^{\bar{D}_f}\in\mathrm{Im}\,\Theta$ is expressed as a polynomial of $f_j\otimes t^{\bar{D}_{f_j}}$ 's and $t^{\bar{E}_i}$'s. By evaluating $t=1$, we obtain an expression of $f$ as a polynomial of $f_j$'s such that the values of $\nu_j$'s ($j\in A$) of each monomial are greater than or equal to those of $f$. Replacing $f_j$'s by the expressions as polynomials of elements in $S$ satisfying ($*A\cup\{i\}$), we obtain an expression of $f$ satisfying ($*A\cup\{i\}$).
\qed

\vspace{3mm}

\begin{prop}\label{prop:(0,i)}
Fix $S=S_p\subset\mathbb{C}[V]^{[G,G]}$ and $i\in\{1,\dots,m\}$. Then the equality $\mathrm{min}_i(I)=\mathrm{min}_i(J)$ holds if and only if $S$ satisfies ($*\{i\}$). Moreover, {\bf Step ($\boldsymbol{0,i}$)} ends in finitely many times (i.e. the equality $\mathrm{min}_i(I)=\mathrm{min}_i(J)$ holds for $S_q$ with $q\gg p$).
\end{prop}

{\em Proof.} First assume that $S$ satisfies ($*\{i\}$). Let $h$ be $\deg_i$-homogeneous element of $\mathrm{min}_i(J)$. Then  there is $\tilde{h}\in\mathbb{C}[X_1,\dots,X_k]$ such that $\alpha(\tilde{h})=\alpha(h)$ and $\deg_i(\tilde{h})>\deg_i(h)$ by Lemma \ref{lem:reduction}. Since $h$ is $\deg_i$-homogeneous, $h$ is the minimal part of $h-\tilde{h}\in I$ and therefore $h\in\mathrm{min}_i(I)$. By using the fact $\mathrm{min}_i(J)$ is a $\deg_i$-homogeneous ideal, we conclude that $\mathrm{min}_i(I)=\mathrm{min}_i(J)$.

Conversely, we assume that $\mathrm{min}_i(I)=\mathrm{min}_i(J)$ holds. Then for any $h\in\mathrm{min}_i(J)$, there is $\tilde{h}\in\mathbb{C}[X_1,\dots,X_k]$ such that $h-\tilde{h}\in I$ and $\deg_i(h)<\deg_i(\tilde{h})$. Thus $h$ satisfies ($*\emptyset,i,p$), and $S$ satisfies ($*\{i\}$) by Lemma \ref{lem:reduction}.

To prove the second claim, it is enough to show that $S_q$ satisfies ($*\{i\}$) with $q\gg p$. By the argument in the proof of Lemma \ref{lem:reduction}, we only have to check the condition ($*\{i\}$) for each $f\in\{f_1,\dots,f_s\}$. Similarly to the ``if'' part of the proof of Lemma \ref{lem:reduction}, take $h'\in\mathbb{C}[X_1,\dots,X_k]$ and $h=\mathrm{min}_i(h')\in\mathrm{min}_i(J)$ such that $\alpha(h')=f$. Then one can write $h=h''+\sum_{j=1}^l a_j h_j$ where $h''\in\mathrm{min}_i(I)$, $a_j\in\mathbb{C}[X_1,\dots,X_k]$, and $h_j$'s are ones in {\bf Step ($\boldsymbol{0,i}$)}. Since $h''\in\mathrm{min}_i(I)$, there is $\tilde{h''}\in\mathbb{C}[X_1,\dots,X_k]$ such that $h''-\tilde{h}''\in I$ and $\deg_i(h'')<\deg_i(\tilde{h}'')$. Let $X_{k_1},\dots,X_{k_l}\in\mathbb{C}[X_1,\dots,X_{|S_{p+1}|}]$ be the new variables associated to $S_{p+1}$ corresponding to $h_1,\dots,h_l$ respectively. Then replacing $h''$ and $h_j$'s in $h$ by $\tilde{h}''$ and $X_{k_j}$ respectively increases the value of $\deg_i$. Repeating this process gives $h'\in\mathbb{C}[X_1,\dots,X_{|S_q|}]$ for $q\gg p$ such that $\nu_i(f)=\deg_i(h')$ by induction on $\nu_i(f)-\deg_i(h')$.
\qed

\vspace{3mm}
\begin{prop}\label{prop:(i',i)}
Fix $S=S_p\subset\mathbb{C}[V]^{[G,G]}$ and $1\le i'<i\le m$. Assume that $S$ satisfies ($*\{1,\dots,i'-1,i\}$). Then the equality $\tilde{I}_{i',i}=\tilde{I}'_{i',i}$ holds if and only if $S$ satisfies ($*\{1,\dots,i',i\}$).
\end{prop}

{\em Proof.} First assume that $S$ satisfies ($*\{1,\dots,i',i\}$). Let $h'$ be any element in $\tilde{I}_{i',i}$ which is homogeneous with respect to $\deg_j$ for each $j=1,\dots,i',i$. If $t_{i'}|h'$ or $t_i|h'$, then clearly $h'$ is in $\tilde{I}'_{i',i}$. So we assume otherwise. Set $h=\mathrm{min}_{i}(h'|_{t=1})\in\mathrm{min}_i(I)$. Note that $\mathrm{min}_i(I)=\mathrm{min}_i(J)$ by Proposition \ref{prop:(0,i)}. Since $S$ satisfies ($*\{1,\dots,i',i\}$), there is $\tilde{h}\in\mathbb{C}[X_1,\dots,X_k]$ such that $h-\tilde{h}\in I$, $\deg_{i}(h)<\deg_{i}(\tilde{h})$, and $\deg_j(h)\le\deg_j(\tilde{h})$ for $j=1,\dots,i'$ by Lemma \ref{lem:reduction}. Set $h''=(h-\tilde{h})_{\{1,\dots,i',i\}}\in I_{\{1,\dots,i',i\}}$. By the choice of $\tilde{h}$ and the homogeneity of $h'$ and $h''$, one sees that $h'|_{t_i=0}= h''|_{t_i=0}\prod_{j=1}^{i'}t_j^{l_j}$ in $R_{1,\dots,i',i}$ for some $l_j\ge0$ ($j=1,\dots,i'-1$) and $l_{i'}>0$. Therefore $h'\in \tilde{I}'_{i',i}$.

Conversely, we assume that $\tilde{I}_{i',i}=\tilde{I}'_{i',i}$ holds. To prove that $S$ satisfies ($*\{1,\dots,i',i\}$), it is enough to show that, for any $\deg_i$-homogeneous $h\in\mathrm{min}_i(J)$, there exists $\tilde{h}\in\mathbb{C}[X_1,\dots,X_k]$ such that $h-\tilde{h}\in I$, $\deg_i(\tilde{h})>\deg_i(h)$, and $\deg_j(\tilde{h})\ge\deg_j(h)$ for $j=1,\dots,i'$ by Lemma \ref{lem:reduction}.

Since $S$ satisfies ($*\{1,\dots,i'-1,i\}$), there is $\tilde{h}\in\mathbb{C}[X_1,\dots,X_k]$ such that $h-\tilde{h}\in I$, $\deg_i(\tilde{h})>\deg_i(h)$, and $\deg_j(\tilde{h})\ge\deg_j(h)$ for $j=1,\dots,i'-1$. Let $h'=(h-\tilde{h})_{\{1,\dots,i',i\}}\in I_{\{1,\dots,i',i\}}$. If $h'|_{t_{i'}=t_i=0}\ne0$ in $\mathbb{C}[X_1,\dots,X_k]$, then $\deg_{i'}(\tilde{h})\ge\deg_{i'}(h)$ and thus we are done. So we assume otherwise. Then $h'$ is in $\tilde{I}_{i',i}=\tilde{I}'_{i',i}$ and we can write $h'=h'_1t_{i'}+h'_2t_i$ with $h'_j\in I_{\{1,\dots,i',i\}}$. Therefore one has $h'|_{t_i=0}=h'_1|_{t_i=0} t_{i'}$. By repeating this process, we should reach the case where $h'_1|_{t_{i'}=t_i=0}\ne0$. In this case $\tilde{h}:=(h'_1-h'_1|_{t_i=0})|_{t=1}\in\mathbb{C}[X_1,\dots,X_k]$ satisfies the desired condition.
\qed

\vspace{3mm}

Consider the situation where one has obtained $S_{p+j}$ from $S_{p+j-1}$ for $j=1,\dots,c$ by performing {\bf Step ($\boldsymbol{i',i}$)} ($i'\ne0$) and assume that $S_p\subsetneq S_{p+1}\subsetneq \cdots \subsetneq S_{p+c}$. For $j=1,\dots,c$, define $B_{p+j,i',i}\subset\mathbb{C}[X_1,\dots,X_{S_{p+j-1}}]$ as the set
$$\begin{aligned}\{h\in\mathrm{min}_i(I)\,|\, &h\text{ is }\deg_i\text{-homogeneous and}\\
&\text{does not satisfy either }(*A,i,p+j-1)\text{ or }(*A,i,p+j)\}
\end{aligned}$$
where $I$ is defined with respect to $S_{p+j-1}$ and $A=\{1,\dots,i'\}$.

\begin{lem}\label{lem:degree bound}
Assume that $S_p$ satisfies ($*\{1,\dots,i'-1,i\}$). Then there are integers $m_1<m_2<\dots<m_c$ such that the inequality $\deg_{i'}(h)>m_j$ holds for any $h\in B_{p+j,i',i}$ and $j=1,\dots,c$. 
\end{lem}

{\em Proof.} Take any element $h\in B_{p+1,i',i}$. Since $S_p$ satisfies ($*A\{1,\dots,i'-1,i\}$), there is $\tilde{h}\in\mathbb{C}[X_1,\dots,X_k]$ such that $h-\tilde{h}\in I$, $\deg_{i}(h)<\deg_{i}(\tilde{h})$, and $\deg_j(h)\le\deg_j(\tilde{h})$ for $j\in \{1,\dots,i'-1,i\}$ by Lemma \ref{lem:reduction}. Then $h':=(h-\tilde{h})_{\{1,\dots,i',i\}}$ is in $\tilde{I}_{i',i}$ since $h$ does not satisfy ($*\{1,\dots,i',i\},i,p$). Therefore we can write $h'=h''+\sum_{j=1}^l a_j h_j$ with $h''\in \tilde{I}'_{i',i}$ and $a_j\in R_{\{1,\dots,i',i\}}$ where $h_j$'s are ones in {\bf Step ($\boldsymbol{i',i}$)}. By the argument of the proof of Proposition \ref{prop:(i',i)}, one sees that there is $\tilde{h}''\in\mathbb{C}[X_1,\dots,X_k]$ such that $\mathrm{min}_i(h''|_{t=1})-\tilde{h}''\in I$, $\deg_{i}(\mathrm{min}_i(h''|_{t=1}))<\deg_{i}(\tilde{h}'')$, and $\deg_j(\mathrm{min}_i(h''|_{t=1}))\le\deg_j(\tilde{h}'')$ for $j\in \{1,\dots,i',i\}$. On the other hand, let $X_{k_1},\dots,X_{k_l}\in\mathbb{C}[X_1,\dots,X_{|S_{p+1}|}]$ be the new variables associated to $S_{p+1}$ corresponding to $h_1,\dots,h_l$ respectively. Then $\mathrm{min}_i(h_j|_{t=1})-X_{k_j}$ is in $I$, $\deg_{i}(\mathrm{min}_i(h_j|_{t=1}))<\deg_{i}(X_{k_j})$, and $\deg_{j'}(\mathrm{min}_i(h_j|_{t=1}))\le\deg_{j'}(X_{k_j})$ for $j'\in\{1,\dots,i'\}$. By the choice of $\tilde{h}$ and the homogeneity of $h_j$'s, the inequality $\deg_{j'}(h)\le\underset{j=1,\dots,l}{\mathrm{min}}\{\deg_{j'}(\mathrm{min}_i(a_j h_j|_{t=1}))\}$ for $j'\in \{1,\dots,i'-1\}$ holds.
 If $\deg_{i'}(h)\le m_1:=\underset{j=1,\dots,l}{\mathrm{min}}\{\deg_{i'}(\mathrm{min}_i(h_j|_{t=1}))\}$ also holds, then
$$\deg_i(h)=\deg_i(\mathrm{min}_i((h''+\sum_{j=1}^l a_j h_j)|_{t=1}))<\deg_i(\tilde{h}''+\sum_{j=1}^l a_j|_{t=1} X_{k_j})$$
and
$$\deg_{i'}(h)=\deg_{i'}(\mathrm{min}_i((h''+\sum_{j=1}^l a_j h_j)|_{t=1}))\le\deg_{i'}(\tilde{h}''+\sum_{j=1}^l a_j|_{t=1} X_{k_j}).$$
This is contrary to the fact that $h$ does not satisfy ($*\{1,\dots,i'\},i,p+1$) and hence $\deg_{i'}(h)>m_1$.

Since $S_{p+1}\subsetneq S_{p+2}$, there are $h_{l+1},\dots,h_{l'}\in\mathbb{C}[X_1,\dots,X_{|S_{p+1}|}]$ such that $\tilde{I}_{i',i}=\tilde{I}'_{i',i}+(h_{l+1},\dots,h_{l'})$ as in {\bf Step ($\boldsymbol{i',i}$)} for $S=S_{p+1}$. By construction $h_{l+1},\dots,h_{l'}$ come from elements in $\tilde{I}_{i',i}\backslash\tilde{I}'_{i',i}$. Therefore they do not satisfy ($*\{1,\dots,i'\},i,p+1$). The same argument as above shows that $m_2:=\underset{j=l+1,\dots,l'}{\mathrm{min}}\{\deg_{i'}(\mathrm{min}_i(h_j|_{t=1}))\}>m_1$. The integers $m_3,\dots,m_c$ are defined similarly and the claim about the inequality for $j=2,\dots,c$ follows from the same argument as the case $j=1$.
\qed

\vspace{3mm}

\begin{prop}\label{prop:termination}
{\bf Step ($\boldsymbol{i',i}$)} ends in finitely many times (i.e. the equality $\tilde{I}_{i',i}=\tilde{I}'_{i',i}$ holds for $S_q$ with $q\gg p$).
\end{prop}
{\em Proof.} 
It is enough to show that $S_q$ satisfies ($*\{1,\dots,i',i\}$) for $q\gg p$ by Proposition \ref{prop:(i',i)}. By the argument in the proof of Lemma \ref{lem:reduction}, we only have to check the condition ($*\{1,\dots,i',i\}$) for each $f\in\{f_1,\dots,f_s\}$. Similarly to the ``if'' part of the proof of Lemma \ref{lem:reduction}, take $h'\in\mathbb{C}[X_1,\dots,X_k]$ and $h=\mathrm{min}_i(h')\in\mathrm{min}_i(J)$ such that $\alpha(h')=f$. By Lemma \ref{lem:degree bound}, $h$ is not in $B_q$ for $q\gg p$. Therefore $h$ satisfies ($*\{1,\dots,i'\},i,q$) for $q\gg p$.
\qed

\begin{cor}\label{cor:finite}
The algorithm ends in finitely many times, that is, $S_{\infty}:=\bigcup_p S_p$ is a finite set. Moreover, the subset $\{\phi\otimes t^{\bar{D}_{\phi}}\}_{\phi\in S_\infty}\cup\{t^{\bar{E}_1},\dots,t^{\bar{E}_m}\}$
of $\mathbb{C}[V]^{[G,G]}\otimes_{\mathbb{C}}\mathbb{C}[\mathrm{Cl}(X)^{\mathrm{free}}]$ is a generating system of $\mathrm{Im}\,\Theta$.
\end{cor}

{\em Proof.} The first claim follows from Proposition \ref{prop:(0,i)} and \ref{prop:termination}. As we see that $S_\infty$ satisfies ($*\{1,\dots,m\}$)(=($*$)) by Proposition \ref{prop:(i',i)}, the second claim follows from Proposition \ref{prop:generator}.
\qed

\vspace{3mm}

This theorem makes it possible for us to calculate Cox rings of minimal models of any quotient singularities at least theoretically. Before we try concrete examples, we state one application of the above construction of Cox rings.

In the previous section, we showed that the simply-connectedness of the regular part of a minimal model is determined by whether $G$ is generated by junior elements. By using the embedding of the Cox ring above, we can show that the torsion freeness of the divisor class group of the minimal model can also be read from the property of $G$.

\begin{prop}
The divisor class group $\mathrm{Cl}(X)$ of the minimal model $X$ of $V/G$ is torsion-free if and only if $G$ is generated by $[G,G]$ and junior elements.
\end{prop}

{\em Proof.} Let $H$ be the subgroup of $G$ generated by $[G,G]$ and junior elements. First assume $H\ne G$. Then there is an element $f\in\mathbb{C}[V]^{H}\backslash\mathbb{C}[V]^{G}\subset\mathbb{C}[V]^{[G,G]}$ which is homogeneous with respect to $Ab(G)$-action. Let $\bar{D}_f\in\mathrm{Cl}(X)^{\mathrm{free}}$ be the divisor class associated to $f$. Then, by Lemma \ref{lem:strict transform}, one has $\bar{D}_f=-\sum_{i=1}^m\frac{1}{r_i} \nu_i(f)\bar{E}_i\in\bigoplus_{i=1}^r \mathbb{Z}\bar{E}_i$. This implies that the integral Weil divisor $D_f+\sum_{i=1}^m\frac{1}{r_i} \nu_i(f)E_i$ is torsion but nonzero in $\mathrm{Cl}(X)$.

Next assume that $H=G$. Let $D$ be any Weil divisor which is a torsion in $\mathrm{Cl}(X)$. Let $f\in\mathbb{C}[V]^{[G,G]}\cong\mathrm{Cox}(V/G)$ be the defining section of $\pi_*D$ and let $D'$ be the strict transform of $\pi_*D$ on $X$. Clearly one can write $D-D'=\sum_{i=1}^m a_i E_i$ for some $a_i\in\mathbb{Z}$. By Lemma \ref{lem:strict transform}, we have $\bar{D}'=-\sum_{i=1}^m\frac{1}{r_i} \nu_i(f)\bar{E}_i$ in $\mathrm{Cl}(X)^{\mathrm{free}}$. On the other hand, we have $\bar{D}'=\bar{D}-\sum_i a_i\bar{ E}_i=-\sum_i a_i \bar{E}_i$ in $\mathrm{Cl}(X)^{\mathrm{free}}$ since $D$ is a torsion. Therefore $-\sum_{i=1}^m\frac{1}{r_i} \nu_i(f)\bar{E}_i=-\sum_i a_i \bar{E}_i$. As $\{E_i\}_i$ is a $\mathbb{Q}$-basis of $\mathrm{Cl}(X)_{\mathbb{Q}}$, the condition $\nu_i(f)\equiv 0\mbox{ mod }r_i\;(\forall i)$ must be satisfied. This is equivalent to $f\in\mathbb{C}[V]^H=\mathbb{C}[V]^G=\mathrm{Cox}(V/G)_0$. Hence the class of $D$ in $\mathrm{Cl}(X)$ is contained in $\bigoplus_{i=1}^r \mathbb{Z}E_i$. However, $\bigoplus_{i=1}^r \mathbb{Z}E_i$ is torsion-free and thus $D$ must be 0.
\qed

\vspace{3mm}
Now we calculate Cox rings for several examples. For this, consider the Laurent polynomial ring $R:=\mathbb{C}[V]^{[G,G]}[t_1^{\pm 1},\dots,t_m^{\pm 1}]$ over $\mathbb{C}[V]^{[G,G]}$. By Proposition \ref{prop:generator}, we can regard $\mathrm{Im}\,\Theta$ as a subring of $R$ by identifying $f\otimes t^{\bar{D}_f}\in\mathrm{Im}\,\Theta$ with 
$ft_1^{\nu_1(f)}\cdots t_m^{\nu_m(f)}$. Most of the calculations in the examples below are done with a computer making use of the softwares ``Macaulay 2''\cite{GS} or ``SINGULAR''\cite{GPS}. See Appendix \ref{App1} for how to perform the algorithm above and see Appendix \ref{App2} for how to calculate the relations of the generators of the Cox rings efficiently.

\vspace{3mm}
{\bf Remark.}\;In the following examples, a finite group $G$ is realized as a subgroup of the matrix group $SL_n(\mathbb{C})(=SL(V))$. The letters $x,y,\dots$ denote the dual basis to the standard basis of $V=\mathbb{C}^n$. Originally, we should let $G$ act on $V^*$ as the dual representation of $G$ on $V$. However, for convenience, we will let $G$ act on $V^*$ by identifying the dual basis of $V^*$ with the standard basis of $V$. This difference will not produce any effect on the result since one representation and its dual give rise to isomorphic quotient singularities.

\vspace{5mm}
\noindent{\bf Example 1.} (ADE-singularities)\\
{\bf Case 1.} $A_m$-singularity ($m\ge 1$)\\
$\bullet\;$ Cyclic group $G=
\langle g_1=\begin{pmatrix}\zeta&0\\0&\zeta^{-1}\end{pmatrix}\rangle,\;\zeta=\mathrm{exp}\,\frac{2\pi\sqrt{-1}}{m+1}$\\
$\bullet$ representatives of junior elements: $g_k:=g_1^k\;(k=1,\dots,m)$.\\
$\bullet$ $x$ (resp. $y$) is a $\zeta$ (resp. $\zeta^{-1}$)-eigenvector of $g_1$.

The order $r_k$ of $g_k$ is given by $r_k=\frac{m+1}{gcd(k,m+1)}$, and the valuations are given by
$$\nu_k(x)=\frac{kr_k}{m+1}\mbox{ and }\nu_k(y)=r_k-\frac{kr_k}{m+1}.$$
Since $xt_1^{\frac{r_1}{m+1}}\cdots t_m^{\frac{mr_m}{m+1}},yt_1^{r_1-\frac{r_1}{m+1}}\cdots t_m^{r_m-\frac{mr_m}{m+1}},t_1^{-r_1},\dots,t_m^{-r_m}\in R$ clearly have no relations, the algorithm trivially means that these are the generators of the Cox ring of the minimal model (or the crepant resolution) of $V/G$.

\vspace{3mm}
\noindent{\bf Case 2.} $D_m$-singularity ($m\ge 4$)\\
$\bullet$ Binary dihedral group $G=
\langle g_1=\begin{pmatrix}\zeta&0\\0&\zeta^{-1}\end{pmatrix},
g_{m-1}=\begin{pmatrix}0&1\\-1&0\end{pmatrix}\rangle,\;
\zeta=\mathrm{exp}\,\frac{2\pi\sqrt{-1}}{2(m-2)}$,\\
\hspace{3mm}$[G,G]=\langle g_1^2\rangle,\;Ab(G)\cong\begin{cases}(\mathbb{Z}/2\mathbb{Z})^{\times 2}&\mbox{if }m\mbox{ is even} \\
\mathbb{Z}/4\mathbb{Z}&\mbox{if }m\mbox{ is odd}\end{cases}$.\\
$\bullet$ representatives of junior elements: $g_k:=g_1^k\;(k=1,\dots,m-2),g_{m-1},g_m:=g_{m-1}g_1$.\\
$\bullet$ $x$ (resp. $y$) is a $\zeta$ (resp. $\zeta^{-1}$)-eigenvector of $g_k\;(k=1,\dots,m-2)$.\\
$\bullet$ $x+\sqrt{-1}y$ (resp. $x-\sqrt{-1}y$) is a $\sqrt{-1}$ (resp. $-\sqrt{-1}$)-eigenvector of $g_{m-1}$.\\
$\bullet$ $x+\sqrt{-1}\zeta y$ (resp. $x-\sqrt{-1}\zeta y$) is a $\sqrt{-1}$ (resp. $-\sqrt{-1}$)-eigenvector of $g_m$.

The order $r_k$ of $g_k$ is given by
$r_k=\begin{cases}\frac{2(m-2)}{gcd(k,2(m-2))}&\mbox{if}\;\;k=1,\dots,m-2 \\4&\mbox{if}\;\;k=m-1,m.\end{cases}$

In this case we can take $x^{m-2}+(\sqrt{-1}y)^{m-2},x^{m-2}-(\sqrt{-1}y)^{m-2}$ and $xy$ as homogeneous generators of $\mathbb{C}[V]^{[G,G]}$. We can directly calculate the valuations as follows.
$$\nu_k(x^{m-2}+(\sqrt{-1}y)^{m-2})=\begin{cases}\frac{kr_k}{2}&\mbox{if}\;\;k=1,\dots,m-2\\m-2&\mbox{if}\;\;k=m-1\\m&\mbox{if}\;\;k=m\end{cases}$$
$$\nu_k(x^{m-2}-(\sqrt{-1}y)^{m-2})=\begin{cases}\frac{kr_k}{2}&\mbox{if}\;\;k=1,\dots,m-2\\m&\mbox{if}\;\;k=m-1\\m-2&\mbox{if}\;\;k=m\end{cases}$$
$$\nu_k(xy)=\begin{cases}r_k&\mbox{if}\;\;k=1,\dots,m-2\\2&\mbox{if}\;\;k=m-1,m\end{cases}$$

Now we apply the algorithm to $S=\{x^{m-2}+(\sqrt{-1}y)^{m-2},x^{m-2}-(\sqrt{-1}y)^{m-2},xy\}$. We use the notations in the algorithm above. First, the kernel of $\alpha:\mathbb{C}[X_1,X_2,X_3]\to\mathbb{C}[V]^{[G,G]}$ is a principal ideal $I=(X_1^2-X_2^2-4(\sqrt{-1}X_3)^{m-2})$ and one sees that
$$\mathrm{min}_k(I)=\begin{cases}(X_1^2-X_2^2)&\mbox{if}\;\;k=1,\dots,m-3\\
(X_1^2-X_2^2-4(\sqrt{-1}X_3)^{m-2})&\mbox{if}\;\;k=m-2\\
(X_1^2-4(\sqrt{-1}X_3)^{m-2})&\mbox{if}\;\;k=m-1\\
(-X_2^2-4(\sqrt{-1}X_3)^{m-2})&\mbox{if}\;\;k=m.\end{cases}$$

On the other hand, one sees that
$$\mathrm{min}_k(x^{m-2}+(\sqrt{-1}y)^{m-2})=\begin{cases}x^{m-2}&\mbox{if}\;\;k=1,\dots,m-3\\
x^{m-2}+(\sqrt{-1}y)^{m-2}&\mbox{if}\;\;k=m-2\\
\frac{1}{2^{m-1}}(x+\sqrt{-1}y)^{m-2}&\mbox{if}\;\;k=m-1\\
\frac{m-2}{2^{m-1}}(x+\sqrt{-1}\zeta y)^{m-3}(x-\sqrt{-1}\zeta y)&\mbox{if}\;\;k=m,\end{cases}$$

$$\mathrm{min}_k(x^{m-2}-(\sqrt{-1}y)^{m-2})=\begin{cases}x^{m-2}&\mbox{if}\;\;k=1,\dots,m-3\\
x^{m-2}-(\sqrt{-1}y)^{m-2}&\mbox{if}\;\;k=m-2\\
\frac{m-2}{2^{m-1}}(x+\sqrt{-1}y)^{m-3}(x-\sqrt{-1}y)&\mbox{if}\;\;k=m-1\\
\frac{1}{2^{m-1}}(x+\sqrt{-1}\zeta y)^{m-2}&\mbox{if}\;\;k=m,\end{cases}$$
and
$$\mathrm{min}_k(xy)=\begin{cases}xy&\mbox{if}\;\;k=1,\dots,m-2\\
\frac{1}{4\sqrt{-1}}(x+\sqrt{-1}y)^2&\mbox{if}\;\;k=m-1\\
\frac{\zeta^{-1}}{4\sqrt{-1}}(x+\sqrt{-1}\zeta y)^2&\mbox{if}\;\;k=m.\end{cases}$$

From these, one can check that $\mathrm{min}_k(I)=\mathrm{min}_k(J)$ for all $k$. One can also check that each {\bf Step ($\boldsymbol{i',i}$)} ends at one try. Therefore, by Corollary \ref{cor:finite}, we obtain a generating system of the Cox ring
$$\begin{aligned}&(x^{m-2}+(\sqrt{-1}y)^{m-2})t_1^{m-2}\cdots t_{m-1}^{m-2}t_{m}^m,(x^{m-2}-(\sqrt{-1}y)^{m-2})t_1^{m-2}\cdots t_{m-2}^{m-2}t_{m-1}^mt_{m}^{m-2},\\
&(xy)t_1^{r_1}\cdots t_{m-2}^{r_{m-2}}t_{m-1}^2t_m^2,t_1^{-r_1},\dots,t_m^{-r_m}.\end{aligned}$$

If we rename these elements as $X_1,X_2,X_3,Y_1,\dots,Y_m$ in order, then they have a single relation $$X_1^2Y_m-X_2^2Y_{m-1}-4(\sqrt{-1}X_3)^{m-2}\prod_{k=1}^{m-3}Y_k^{m-2-k}=0.$$

\vspace{3mm}
\noindent{\bf Case 3.} $E_6$-singularity \\
$\bullet$ Binary tetrahedral group $G=
\langle g_1=\begin{pmatrix}\sqrt{-1}&0\\0&\sqrt{-1}\end{pmatrix},
g_2=\frac{1}{\sqrt{2}}\begin{pmatrix}\zeta&\zeta\\\zeta^3&\zeta^7\end{pmatrix}\rangle,\\
\hspace{3mm}\zeta=\mathrm{exp}\,\frac{2\pi\sqrt{-1}}{8},\;[G,G]=\langle g_1,\begin{pmatrix}0&1\\-1&0\end{pmatrix}\rangle,\;
Ab(G)\cong\mathbb{Z}/3\mathbb{Z}$.\\
$\bullet$ representatives of junior elements: $g_1,g_k:=g_2^{k-1}\;(k=2,\dots,6)$.\\
$\bullet$ $x$ (resp. $y$) is a $\sqrt{-1}$ (resp. $\sqrt{-1}$)-eigenvector of $g_1$.\\
$\bullet$ $x+(\sqrt{2}\omega\zeta^3-1)y$ (resp. $x+(\sqrt{2}\omega\zeta^7-\zeta^2)y$) is a $-\omega^{k-1}$(resp. $-\omega^{-k+1}$)-eigenvector\\
\hspace{3mm}of $g_k,\;(k=2,\dots,6)$ where $\omega=\mathrm{exp}\,\frac{2\pi\sqrt{-1}}{3}$.\\

We can take $x^4+y^4+2\sqrt{-3}x^2y^2,\,x^4+y^4+2\sqrt{-3}x^2y^2$ and $x^5y-xy^5$ as homogeneous generators of $\mathbb{C}[V]^{[G,G]}$. The information of the valuations are summarized as follows.

\begin{center}
\begin{tabular}{|c|c|c|c|c|c|c|}\hline
$k$&1&2&3&4&5&6\\ \hline \hline
$\nu_k(x^4+y^4+2\sqrt{-3}x^2y^2)$&4&4&4&4&5&8\\ \hline 
$\nu_k(x^4+y^4-2\sqrt{-3}x^2y^2)$&4&8&5&4&4&4\\ \hline
$\nu_k(x^5y-xy^5)$&8&6&6&6&6&6\\ \hline
$r_k$&4&6&3&2&3&6\\ \hline
\end{tabular}
\end{center}

By applying the algorithm to $S=\{x^4+y^4+2\sqrt{-3}x^2y^2,\,x^4+y^4-2\sqrt{-3}x^2y^2,\,x^5y-xy^5\}$, one sees that $$\begin{aligned}
&(x^4+y^4+2\sqrt{-3}x^2y^2)t_1^4t_2^8t_3^5t_4^4t_5^4t_6^4,
(x^4+y^4-2\sqrt{-3}x^2y^2)t_1^4t_2^4t_3^4t_4^4t_5^5t_6^8,\\
&(x^5y-xy^5)t_1^8t_2^6t_3^6t_4^6t_5^6t_6^6,
t_1^{-4},t_2^{-6},t_3^{-3},t_4^{-2},t_5^{-3},t_6^{-6}
\end{aligned}$$
in $R$ are generators of the Cox ring.

If we rename these elements as $X_1,X_2,X_3,Y_1,\dots,Y_6$ in order, then they have a single relation $$X_1^3Y_2^2Y_3-X_2^3Y_5Y_6^2-12\sqrt{-3}X_3^2Y_1=0.$$

\vspace{3mm}
\noindent{\bf Case 4.} $E_7$-singularity \\
$\bullet$ Binary octahedral group $G=
\langle g_1=\begin{pmatrix}\zeta&0\\0&\zeta^{-1}\end{pmatrix},
g_5=\frac{1}{\sqrt{2}}
\begin{pmatrix}\zeta&\zeta\\\zeta^3&\zeta^7
\end{pmatrix}\rangle,\;\zeta=\mathrm{exp}\,\frac{2\pi\sqrt{-1}}{8}$,\\
\hspace{3mm}$[G,G]=\langle g_1^2,g_2\rangle$(=binary tetrahedral group),\;
$Ab(G)\cong\mathbb{Z}/2\mathbb{Z}$.\\
$\bullet$ representatives of junior elements: $g_1,g_2:=g_1^2,g_3:=g_1^3,g_4:=g_1^4,g_5,g_6:=g_5^2$,\\
\hspace{3mm}$g_7=
\frac{1}{\sqrt{2}}\begin{pmatrix}\sqrt{-1}&1\\-1&-\sqrt{-1}\end{pmatrix}$.\\
$\bullet$ $x$ (resp. $y$) is a $\sqrt{-1}$ (resp. $\sqrt{-1}$)-eigenvector$\sqrt{-1}$-eigenvector of $g_1$.\\
$\bullet$ $x+(\sqrt{2}\omega\zeta^3-1)y$ and $x+(\sqrt{2}\omega\zeta^7-\zeta^2)y$ are eigenvectors of $g_5$ and $g_6$ where $\omega=\\
\hspace{3mm}\mathrm{exp}\,\frac{2\pi\sqrt{-1}}{3}$.\\
$\bullet$ $x-(1-\sqrt{2})\sqrt{-1}y$ (resp. $x-(1+\sqrt{2})\sqrt{-1}y$) is a $\sqrt{-1}$ (resp. $-\sqrt{-1}$)-eigenvector of\\
\hspace{3mm}$g_7$.\\

We can take $x^{12}-33x^8y^4-33x^4y^8+y^{12},\,x^8+14x^4y^4+y^8$ and $x^5y-xy^5$ as homogeneous generators of $\mathbb{C}[V]^{[G,G]}$. The information of the valuations are summarized as follows.

\begin{center}
\begin{tabular}{|c|c|c|c|c|c|c|c|}\hline
$k$&1&2&3&4&5&6&7\\ \hline \hline
$\nu_k(x^{12}-33x^8y^4-33x^4y^8+y^{12})$&12&12&36&12&12&12&14\\ \hline 
$\nu_k(x^8+14x^4y^4+y^8)$&8&8&24&8&12&9&8\\ \hline
$\nu_k(x^5y-xy^5)$&12&8&20&6&6&6&6\\ \hline
$r_k$&8&4&8&2&6&3&4\\ \hline
\end{tabular}
\end{center}

By applying the algorithm to $S=\{x^{12}-33x^8y^4-33x^4y^8+y^{12},\,x^8+14x^4y^4+y^8,\,x^5y-xy^5\}$, one sees that $$\begin{aligned}
&(x^{12}-33x^8y^4-33x^4y^8+y^{12})
t_1^{12}t_2^{12}t_3^{36}t_4^{12}t_5^{12}t_6^{12}t_7^{14},
(x^8+14x^4y^4+y^8)t_1^8t_2^8t_3^{24}t_4^8t_5^{12}t_6^9t_7^8,\\
&(x^5y-xy^5)t_1^{12}t_2^8t_3^{20}t_4^6t_5^6t_6^6t_7^6,
t_1^{-8},t_2^{-4},t_3^{-8},t_4^{-2},t_5^{-6},t_6^{-3},t_7^{-4}
\end{aligned}$$
in $R$ are generators of the Cox ring.

If we rename these elements as $X_1,X_2,X_3,Y_1,\dots,Y_7$ in order, then they have a single relation $$X_1^2Y_7-X_2^3Y_5^2Y_6-108X_3^4Y_1^3Y_2^2Y_3=0.$$

\vspace{3mm}
\noindent{\bf Case 5.} $E_8$-singularity \\
$\bullet$ Binary icosahedral group\\
\hspace{3mm}$G=\langle g_1=\begin{pmatrix}
\epsilon&0\\0&\epsilon^{-1}\end{pmatrix},
h=\frac{1}{\sqrt{5}}\begin{pmatrix}
\hspace{3mm}\zeta^2-\zeta^8&\zeta^4-\zeta^6\\
\zeta^4-\zeta^6&\zeta^8-\zeta^2\end{pmatrix},
g_8=\begin{pmatrix}
0&1\\-1&0\end{pmatrix}\rangle,\;\epsilon=\mathrm{exp}\,\frac{2\pi\sqrt{-1}}{10},\\
\hspace{3mm}\zeta=\mathrm{exp}\,\frac{2\pi\sqrt{-1}}{8},[G,G]=G$.\\
$\bullet$ representatives of junior elements: $g_1,g_2:=g_1^2,g_3:=g_1^3,g_4:=g_1^4,g_5:=g_1^5,g_6:=\\
\hspace{3mm}-g_1h,g_7:=g_6^2,g_8$.\\
$\bullet$ $x$ and $y$ are eigenvectors of $g_1,\dots,g_5$.\\
$\bullet$ $x+(\zeta^3-\omega\zeta^2-\zeta^2-\omega\zeta-1)y$ (resp.  $x+(\zeta^3+\omega\zeta^2+\zeta^2+\omega\zeta-1)y$) are eigenvectors of $g_5$ and $g_6$ where $\omega=\mathrm{exp}\,\frac{2\pi\sqrt{-1}}{3}$.\\
$\bullet$ $x\pm\sqrt{-1}y$ are eigenvectors of $g_8$.\\ 

We can take $x^{30}+y^{30}+522(x^{25}y^5-x^5y^{25})-10005(x^{20}y^{10}+x^{10}y^{20}),\,x^{20}+y^{20}-228(x^{15}y^5-x^5y^{15})+494x^{10}y^{10}$ and $xy(x^{10}+11x^5y^5-y^{10})$ as homogeneous generators of $\mathbb{C}[V]^{[G,G]}$. The information of the valuations are summarized as follows.

\begin{center}
\begin{tabular}{|c|c|c|c|c|c|c|c|c|}\hline
$k$&1&2&3&4&5&6&7&8\\ \hline \hline
$\begin{aligned}\nu_k(&x^{30}+y^{30}+522(x^{25}y^5-x^5y^{25})\\
&-10005(x^{20}y^{10}+x^{10}y^{20}))\end{aligned}$
&30&30&90&60&30&30&30&32\\ \hline 
$\nu_k(x^{20}+y^{20}-228(x^{15}y^5-x^5y^{15})+494x^{10}y^{10})$
&20&20&60&40&20&24&21&20\\ \hline
$\nu_k(xy(x^{10}+11x^5y^5-y^{10}))$
&20&15&40&25&12&12&12&12\\ \hline
$r_k$&10&5&10&5&2&6&3&4\\ \hline
\end{tabular}
\end{center}

By applying the algorithm to $S=\{x^{30}+y^{30}+522(x^{25}y^5-x^5y^{25})-10005(x^{20}y^{10}+x^{10}y^{20}),\,x^{20}+y^{20}-228(x^{15}y^5-x^5y^{15})+494x^{10}y^{10},\,xy(x^{10}+11x^5y^5-y^{10})\}$, one sees that
$$\begin{aligned}
&(x^{30}+y^{30}+522(x^{25}y^5-x^5y^{25})-10005(x^{20}y^{10}+x^{10}y^{20}))
t_1^{30}t_2^{30}t_3^{90}t_4^{60}t_5^{30}t_6^{30}t_7^{30}t_8^{32},\\
&(x^{20}+y^{20}-228(x^{15}y^5-x^5y^{15})+494x^{10}y^{10})
t_1^{20}t_2^{20}t_3^{60}t_4^{40}t_5^{20}t_6^{24}t_7^{21}t_8^{20},\\
&xy(x^{10}+11x^5y^5-y^{10})
t_1^{20}t_2^{15}t_3^{40}t_4^{25}t_5^{12}t_6^{12}t_7^{12}t_8^{12},
t_1^{-10},t_2^{-5},t_3^{-10},t_4^{-5},t_5^{-2},t_6^{-6},t_7^{-3},t_8^{-4}
\end{aligned}$$
in $R$ are generators of the Cox ring.

If we rename these elements as $X_1,X_2,X_3,Y_1,\dots,Y_8$ in order, then they have a single relation $$X_1^2Y_8-X_2^3Y_6^2Y_7-1728X_3^5Y_1^4Y_2^3Y_3^2Y_4=0.$$

\vspace{3mm}
{\bf Remark} By taking linear changes of coordinates, one can check that these results agree with the results in \cite{FGAL} and \cite{D}.

\vspace{3mm}
\noindent{\bf Example 2.} (a group of order 32 acting on a 4-dimensional vector space cf. \cite{BS1},\cite{DW})\\
$\bullet\;G=\langle g_1=\begin{pmatrix}
1&0&0&0\\0&-1&0&0\\0&0&1&0\\0&0&0&-1\end{pmatrix},
g_2=\begin{pmatrix}
0&\sqrt{-1}&0&0\\-\sqrt{-1}
&0&0&0\\0&0&0&-\sqrt{-1}\\0&0&\sqrt{-1}&0\end{pmatrix},\\
\hspace{3mm}g_3=\begin{pmatrix}
0&1&0&0\\1&0&0&0\\
0&0&0&1\\0&0&1&0\end{pmatrix},
g_4=\begin{pmatrix}
0&0&0&1\\0&0&-1&0\\
0&-1&0&0\\1&0&0&0\end{pmatrix}\rangle,[G,G]=\langle-\mathrm{Id}\rangle,\;Ab(G)\cong(\mathbb{Z}/2\mathbb{Z})^{\times4}$\\
$\bullet$ representatives of junior elements: $g_1,g_2,g_3,g_4,g_5:=g_1g_2g_3g_4$\\

 We can take $\phi_{12}=-2(xw+yz),\,\phi_{13}=2\sqrt{-1}(-xw+yz),\,\phi_{14}=2\sqrt{-1}(xy+zw),\,\phi_{15}=2(-xy+zw),\,\phi_{23}=2(xz-yw),\,\phi_{24}=-x^2-y^2+z^2+w^2,\,\phi_{25}=\sqrt{-1}(x^2+y^2+z^2+w^2),\,\phi_{34}=\sqrt{-1}(-x^2+y^2-z^2+w^2),\,\phi_{35}=x^2-y^2-z^2+w^2$ and $\phi_{45}=2(xz+yw)$ as homogeneous generators of $\mathbb{C}[V]^{[G,G]}$. The information of the valuations are summarized as follows (cf. \cite[3.13]{DW}).

\begin{center}
\begin{tabular}{|c|c|c|c|c|c|}\hline
$k$&1&2&3&4&5\\ \hline \hline
$\nu_k(\phi_{12})$&1&1&0&0&0\\ \hline 
$\nu_k(\phi_{13})$&1&0&1&0&0\\ \hline
$\nu_k(\phi_{14})$&1&0&0&1&0\\ \hline
$\nu_k(\phi_{15})$&1&0&0&0&1\\ \hline
$\nu_k(\phi_{23})$&0&1&1&0&0\\ \hline 
$\nu_k(\phi_{24})$&0&1&0&1&0\\ \hline
$\nu_k(\phi_{25})$&0&1&0&0&1\\ \hline
$\nu_k(\phi_{34})$&0&0&1&1&0\\ \hline 
$\nu_k(\phi_{35})$&0&0&1&0&1\\ \hline
$\nu_k(\phi_{45})$&0&0&0&1&1\\ \hline
$r_k$&2&2&2&2&2\\ \hline
\end{tabular}
\end{center}

By applying the algorithm to $S=\{\phi_{12},\dots,\phi_{45}\}$, one sees that each step ends at one try and thus
$$\{\phi_{i,j}t_it_j\}_{1\le i<j\le5}\cup\{t_i^{-2}\}_{i=1,\dots,5}$$
in $R$ are generators of the Cox ring as stated in \cite{DW}. The relations of these elements are calculated in \cite{DW}.

\vspace{3mm}
\noindent{\bf Example 3.} (The complex reflection group $G_4$ cf. \cite{Bel},\cite{LS})\\
$\bullet\;G=\langle g_1=-\frac{1}{2}\begin{pmatrix}
(1+\sqrt{-1})\omega&(1+\sqrt{-1})\omega&0&0\\
(-1+\sqrt{-1})\omega&(1-\sqrt{-1})\omega&0&0\\
0&0&(1-\sqrt{-1})\omega^2&(1-\sqrt{-1})\omega^2\\
0&0&(-1-\sqrt{-1})\omega^2&(1+\sqrt{-1})\omega^2
\end{pmatrix},\\
\;g_2=-\frac{1}{2}\begin{pmatrix}
(1+\sqrt{-1})\omega&(1-\sqrt{-1})\omega&0&0\\
(-1-\sqrt{-1})\omega&(1-\sqrt{-1})\omega&0&0\\
0&0&(1-\sqrt{-1})\omega^2&(1+\sqrt{-1})\omega^2\\
0&0&(-1-\sqrt{-1})\omega^2&(-1+\sqrt{-1})\omega^2
\end{pmatrix}\rangle,\\
\;\omega=\mathrm{exp}\,\frac{2\pi\sqrt{-1}}{3},\\
\;[G,G]=\langle
\begin{pmatrix}\sqrt{-1}&0&0&0\\0&-\sqrt{-1}&0&0\\
0&0&-\sqrt{-1}&0\\0&0&0&-\sqrt{-1}
\end{pmatrix},
\begin{pmatrix}0&1&0&0\\-1&0&0&0\\0&0&0&1\\0&0&-1&0\end{pmatrix}\rangle$:the quaternion group\\
$\bullet$ representatives of junior elements: $g_1,g_2$

The homogeneous generators of the invariant ring $\mathbb{C}[x,y,z,w]^{[G,G]}$ are listed as follows:
$\phi_1 = xz+yw,\\
\phi_2 = x^5y-xy^5,\\
\phi_3 = z^5w-zw^5,\\
\phi_4 = x^4+(4a-2)x^2y^2+y^4,\\
\phi_5 = z^4-(4a-2)z^2w^2+w^4,\\
\phi_6 = xw^3+(-2a+1)yzw^2+(2a-1)xz^2w-yz^3,\\
\phi_7 = x^3w-(2a-1)xy^2w+(2a-1)x^2yz-y^3z,\\
\phi_8 = x^2yw^3-x^3zw^2-y^3z^2w+xy^2*z^3,\\ \phi_9=3ax^2w^2-(a-2)x^2z^2+(4a-8)xyzw-(a-2)y^2w^2+3ay^2z^2,  \\
\phi_{10} = z^4+(4a-2)z^2w^2+w^4,\\
\phi_{11} = x^3w+(2a-1)xy^2w-(2a-1)x^2yz-y^3z,\\
\phi_{12} = 5x^4yw-x^5z-y^5w+5xy^4z,\\
\phi_{13}= xyz^4+2x^2zw^3-2y^2z^3w-xyw^4,\\
\phi_{14} = 3(a-1)x^2w^2-(a+1)y^2w^2+4(a+1)xyzw-(a+1)x^2z^2+3(a-1)y^2z^2,\\
\phi_{15} = x^4-(4a-2)x^2y^2+y^4,\\
\phi_{16} = xw^3+(2a-1)yzw^2-(2a-1)xz^2w-yz^3,\\
\phi_{17} = xz^5-5xzw^4-5yz^4w+yw^5,\\
\phi_{18} = 2x^3yw^2-x^4zw+y^4zw-2xy^3z^2$\\
where $a$ denotes $\mathrm{exp}\,\frac{2\pi\sqrt{-1}}{6}$.

The information of the valuations are summarized as follows.

$$\nu_1(\phi_i)=\begin{cases}0&\mbox{if }1\le i\le8\\2&\mbox{if }9\le i\le13\\1&\mbox{if }14\le i\le18\end{cases}$$
$$\nu_2(\phi_i)=\begin{cases}0&\mbox{if }1\le i\le8\\1&\mbox{if }9\le i\le13\\2&\mbox{if }14\le i\le18\end{cases}$$
$$\sharp\langle g_1\rangle=\sharp\langle g_2\rangle=3.$$

We apply the algorithm to $S=\{\phi_1,\dots,\phi_{18}\}$. In this case $\mathrm{min}_1(J)$ is strictly bigger than $\mathrm{min}_1(I)$ and one can check that there is $h=X_1^3+(-6a+3)X_8\in\mathbb{C}[X_1,\dots,X_{18}]$ such that
$$\mathrm{min}_1(J)=\mathrm{min}_1(I)+(h).$$
Thus we should add $$\begin{aligned}\phi_{19}:=\alpha(h)=&x^3z^3+(-6a+3)xy^2z^3+3x^2yz^2w+(6a-3)y^3z^2w+(6a-3)x^3zw^2\\
&+3xy^2zw^2+(-6a+3)x^2yw^3+y^3w^3\end{aligned}$$
to $S$ and try 
{\bf Step ($\boldsymbol{0,1}$)} again.

One can check that each step ends at one try for $S=\{\phi_1,\dots,\phi_{19}\}$ and thus
$$\phi_1,\dots,\phi_8,\phi_9t_1^2t_2,\dots,\phi_{13}t_1^2t_2,\phi_{14}t_1t_2^2,\dots,\phi_{18}t_1t_2^2,\phi_{19}t_1^3t_2^3,t_1^{-3},t_2^{-3}$$
are the generators of the Cox ring.

\vspace{3mm}
{\bf Remark} This result shows that the conjecture in \cite[\S 6]{DG} is negative. One more generator is necessary.

\section{GIT chambers and ample cones\label{5}}
In this section we summarize the basic notions and results about GIT (Geometric Invariant Theory) for the case of the Cox ring of a minimal model $X$ of the quotient singularity $V/G$.

As stated in section \ref{2}, the algebraic torus $T=\mathrm{Hom(\mathrm{Cl}(X),\mathbb{C}^*)}$ acts on the spectrum $\mathfrak{X}=\mathrm{Spec(\mathrm{Cox}(X))}$, and every divisor class $D\in\mathrm{Cl}(X)$ can be considered as a character of $T$. Now we introduce the notion of (semi)stability. Consider the following vector space
$$R(D):=\{f\in H^0(\mathfrak{X},\mathcal{O}_\mathfrak{X})|t\cdot f=D(t)f\mbox{ for all }t\in T\}(=H^0(X,\mathcal{O}_X(D))).$$

\begin{dfn}
We say that a point $x\in\mathfrak{X}$ is $D$-{\bf semistable} if there exist $i\in\mathbb{Z}_{>0}$ and $f\in R(iD)$ such that $f(x)\ne0$. If moreover $x$ has a finite stabilizer and the $T$-orbit of $x$ is closed in $\{x\in\mathfrak{X}|f(x)\ne0\}$, we say that $x$ is $D$-{\bf stable}. $\mathfrak{X}_D^{ss}$ (resp. $\mathfrak{X}_D^s$) denotes the subset of $D$-semistable (resp. $D$-stable) points in $\mathfrak{X}$. We call a divisor class $D$ in $\mathrm{Cl}(X)$ {\bf generic} if $\mathfrak{X}_D^{ss}=\mathfrak{X}_D^s$.
\end{dfn}

We define the {\em GIT-quotient} of $\mathfrak{X}$ by $T$ with respect to $D$ as
$$\mathfrak{X}/\!/_D T:=\mathrm{Proj}\bigoplus_{i=0}^\infty R(iD).$$
Note that there is a natural map from $\mathfrak{X}/\!/_D T$ to $\mathfrak{X}/\!/_0 T=V/G$ for each $D$. GIT quotients of $\mathfrak{X}$ by $T$ have the following property \cite{MFK}.

\begin{prop}\label{prop:git}
The morphism $q:\mathfrak{X}_D^{ss}\to\mathfrak{X}/\!/_D T$ induced by the inclusion $R(iD)\hookrightarrow H^0(\mathcal{O}_\mathfrak{X})$ is a categorical quotient.
Moreover, there is an open subset $U$ of $\mathfrak{X}/\!/_D T$ such that $q^{-1}(U)=\mathfrak{X}_D^s$ and $q|_{\mathfrak{X}_D^s}:\mathfrak{X}_D^s\to U$ is a geometric quotient.
\end{prop}

If two divisor classes in $\mathrm{Cl}(X)$ give the same semistable locus in $\mathfrak{X}$ and hence give the  canonically isomorphic GIT quotients, we call them {\em GIT equivalent}. It is known that GIT equivalence classes give a chamber structure on the finite dimensional real vector space $\mathrm{Cl}(X)_\mathbb{R}:=\mathrm{Cl}(X)\otimes\mathbb{R}$ (cf. \cite[2.3]{T}) i.e.\\
(i) there are only finitely many GIT equivalence classes\\
(ii) for every GIT equivalence class $C$, the closure $\overline{C}$ is a rational polyhedral cone in $\mathrm{Cl}(X)_\mathbb{R}$ and $C$ is a relative interior of $\overline{C}$.

We call $C$ a {\em GIT chamber} if $C$ is not contained in any hyperplane in $\mathrm{Cl}(X)_\mathbb{R}$. The fan on $\mathrm{Cl}(X)_\mathbb{R}$ given by the closures of all GIT chambers is called the {\em GIT fan}. It is known that $D\in\mathrm{Cl}(X)_\mathbb{R}$ is generic if and only if $D$ is in a GIT chamber.

GIT chambers in $\mathrm{Cl}(X)_\mathbb{R}$ are closely related with the birational geometry of $X$. To see this, we introduce ($\pi$-)movable line bundles for the morphism $\pi:X\to V/G$.

\begin{dfn}
A line bundle $L$ on $X$ is $(\boldsymbol{\pi}$-)\textbf{movable} if $\mathrm{codim\,Supp(Coker}\,\alpha)\ge2$ where $\alpha:\pi^*\pi_*L\to L$ is the natural map of sheaves on $X$. The $(\boldsymbol{\pi}$-)\textbf{movable cone} $\mathrm{Mov}(\pi)$ in $\mathrm{Pic}(X)_\mathbb{R}:=\mathrm{Pic}(X)\otimes_{\mathbb{Z}}\mathbb{R}$ is the cone generated by the classes of $\pi$-movable line bundles.
\end{dfn}

Let $\pi':X'\to V/G$ be another minimal model. Then the birational map $\pi'\circ\pi^{-1}:X\dashrightarrow X'$ is an isomorphism in codimension 1 over $V/G$ (see e.g. \cite{Y}; Lemma 3.1). Therefore there is a natural isomorphism $\mathrm{Pic}(X')_\mathbb{R}\to\mathrm{Pic}(X)_\mathbb{R}$, and we call the image of $\mathrm{Amp}(X')$ by this map the {\em ample cone} of $\pi'$. Ample cones satisfy following properties.\\
(1) $\mathrm{Amp}(X')$ and $\mathrm{Amp}(X'')$ are disjoint for different minimal models $X'$ and $X''$.\\
(2) the movable cone is covered with the closures of all ample cones
$$\mathrm{Mov}(\pi)=\bigcup_{X'}\overline{\mathrm{Amp}}(X')$$
where $X'$ runs through all minimal models.

Note that we can regard the movable cone and the ample cones as the subset of $\mathrm{Cl}(X)_\mathbb{R}$ since minimal models are $\mathbb{Q}$-factorial. In \cite{HK}, it is proved that each ample cone $\mathrm{Amp}(X')$ coincides with some GIT chamber and that the corresponding GIT quotient is $X'$. Therefore every minimal model is realized as the GIT quotient $\mathfrak{X}/\!/_D T$ for some divisor $D$ on $X$.

When we are given explicit generators of the Cox ring, it is hard to determine the chamber structure on $\mathrm{Cl}(X)_\mathbb{R}$ directly from the definition. Thus we give a relatively simple way of doing it. To this end, we prepare the following notations. For a finite subset $S\subset\mathrm{Cox}(X)$ consisting of homogeneous elements $\psi_1,\dots,\psi_k$, let $\mathrm{cone}(S)\subset\mathrm{Cl}(X)_\mathbb{R}$ be the cone generated by $\deg(\psi_1),\dots,\deg(\psi_k)\in\mathrm{Cl}(X)_\mathbb{R}$. We also define $\mathcal{F}(S)$ as the coarsest fan on $\mathrm{Cl}(X)_\mathbb{R}$ whose cones refine all the cones of the form $\mathrm{cone}(S')$ for $S'\subset S$. In other words, the union of all proper faces of $\mathcal{F}(S)$ is equal to the union of all proper faces of all cones of the form $\mathrm{cone}(S')$ for $S'\subset S$. Let $\mathcal{Z}(S)\subset\mathfrak{X}$ be the zero locus of the product $\psi_1\cdots\psi_k$. We have the following.

\begin{prop}\label{prop:chamber}
Let $S$ be a generating system of homogeneous elements of $\mathrm{Cox}(X)$ and let $C$ be a chamber in $\mathcal{F}(S)$. Then for any integral divisor $D\in C$, the complement of $D$-semistable locus in $\mathfrak{X}$ is given by
$\bigcap_{S'}\mathcal{Z}(S')$ where $S'$ runs through all subsets of $S$ such that $|S'|=m:=\mathrm{dim}_\mathbb{R}\mathrm{Cl}(X)_\mathbb{R}$ and that $\mathrm{cone}(S')$ contains $D$. In particular, $\mathcal{F}(S)$ is a refinement of the GIT fan.
\end{prop}

{\em Proof.} We first prove the inclusion $\mathfrak{X}\backslash\mathfrak{X}_D^{ss}\subset\bigcap_{S'}\mathcal{Z}(S')$. Let $S'=\{\psi_1,\dots,\psi_m\}$ be any subset of $S$ such that $\mathrm{cone}(S')$ contains $D$. Then there are positive integers $k,a_1,\dots,a_m$ such that $kD=a_1D_1+\cdots+a_mD_m$ where $D_1,\dots,D_m$ are degrees of $\psi_1,\dots,\psi_m$ respectively. If $x\in\mathfrak{X}\backslash\mathfrak{X}_D^{ss}$, by definition $\psi_1^{a_1}\cdots\psi_m^{a_m}(x)=0$ and hence $x$ is in $\mathcal{Z}(S')$.

Conversely, we assume that $x\in\bigcap_{S'}\mathcal{Z}(S')$. Let $k$ be any positive integer and $f$ be any element in $R(kD)$. When we write $f=\sum_j f_j$ where $f_j$'s are monomials of elements in $S$, each $f_j$ must contain a factor which is the product of $m$ elements in $S$ whose degrees are linearly independent since $kD$ is not on any faces of $\mathcal{F}(S)$. Thus $f$ vanishes at $x$.
\qed

\vspace{3mm}
We can also determine the location of the movable cone in the case of minimal models of quotient singularities.

\begin{prop}\label{movable}
Let $\phi_1,\dots,\phi_k\in\mathbb{C}[V]^{[G,G]}$ be elements which are homogeneous with respect to $Ab(G)$-action such that $\{\phi_j\otimes t^{\bar{D}_{\phi_j}}\}_{j=1,\dots,k}\cup\{t^{\bar{E}_1},\dots,t^{\bar{E}_m}\}\subset\mathbb{C}[V]^{[G,G]}\otimes_{\mathbb{C}}\mathbb{C}[\mathrm{Cl}(X)^{\mathrm{free}}]$ is a generating system of the Cox ring of $X$ (cf. Proposition \ref{prop:generator}). Then $\mathrm{cone}(\{\bar{D}_{\phi_j}\}_{j=1,\dots,k})\subset\mathrm{Cl}(X)_\mathbb{R}$ is the movable cone. 
\end{prop}

{\em Proof.} First note that if one takes any divisor $D_0$ from the interior of $\mathrm{cone}(\bar{E}_1,\dots,\bar{E}_m)$, the $D_0$-semistable locus is defined by $\bigcap_{i=1}^m\{t^{\bar{E}_i}\ne0\}\subset\mathfrak{X}$. Since $D_0$ is an effective exceptional divisor on $X$, one can also check by definition that
$$\mathfrak{X}/\!/_{D_0} T=\mathrm{Proj}\bigoplus_{i=0}^\infty H^0(X,\mathcal{O}_X(iD_0))=\mathrm{Proj}\bigoplus_{i=0}^\infty H^0(X,\mathcal{O}_X)=\mathfrak{X}/\!/_0 T=V/G.$$

Let $D\in\mathrm{Cl}(X)_\mathbb{R}$ be a divisor in a GIT chamber. Note that $\mathfrak{X}/\!/_D T$ is a minimal model if and only if it contains $m$ exceptional divisors. (The genericity of $D$ ensures that the GIT quotient is $\mathbb{Q}$-factorial. cf.\cite[1.11(2)]{HK}) Since $D$ gives a geometric quotient, the image of the divisor $\{t^{\bar{E}_i}=0\}\subset\mathfrak{X}$ by the quotient map $\mathfrak{X}^{D-ss}\to\mathfrak{X}/\!/_D T$ is also a divisor as long as $\{t^{\bar{E}_i}=0\}\cap\mathfrak{X}^{D-ss}\ne \emptyset$. In that case the image is nothing but the exceptional divisor on $X$. Therefore it is enough to show that $D\in \mathrm{cone}(\{\bar{D}_{\phi_j}\}_{j=1,\dots,k})$ if and only if $\{t^{\bar{E}_i}=0\}\cap\mathfrak{X}^{D-ss}\ne \emptyset$ for all $i$.

If $D$ is not in $\mathrm{cone}(\{\bar{D}_{\phi_j}\}_{j=1,\dots,k})$, one can check that there is $i$ such that $\mathfrak{X}\backslash\mathfrak{X}^{D-ss}$ contains $\{t^{\bar{E}_i}=0\}$ by Proposition \ref{prop:chamber}. Conversely, if $D$ is in $\mathrm{cone}(\{\bar{D}_{\phi_j}\}_{j=1,\dots,k})$, again by Proposition \ref{prop:chamber}, we see that the non-semistable locus $\mathfrak{X}\backslash\mathfrak{X}^{D-ss}$ is described as a combination of unions and intersections of the zero loci of $\phi_j\otimes t^{\bar{D}_{\phi_j}}$'s. However, the zero locus of $\phi_j\otimes t^{\bar{D}_{\phi_j}}$ does not contain the divisor $\{t^{\bar{E}_i}=0\}$ since the zero locus of $\phi_j\otimes t^{\bar{D}_{\phi_j}}$ on $X$ is the strict transform of a divisor on $V/G$ (see Lemma \ref{lem:strict transform}). Therefore $\mathfrak{X}\backslash\mathfrak{X}^{D-ss}$ can never contain $\{t^{\bar{E}_i}=0\}$.
\qed
\vspace{3mm}

Now we apply the above proposition to a minimal model $X$ of $V/G$ where $G$ is the binary tetrahedral group treated in section \ref{4}, Example 3. According to the results in the previous section, the degrees of the generators of the Cox ring of $X$ are $(0,0),(2,1),(1,2),(3,3),(-3,0)$ and $(0,-3)$. By the proposition, the following figure gives a refinement of the GIT fan on $\mathrm{Cl}(X)_{\mathbb{R}}\cong\mathbb{R}^2$.

\begin{center}
\setlength\unitlength{1truecm}
\begin{picture}(6,6)(0,0)
\put(3.1,3.1){\line(-1,0){3}}
\put(3.1,3.1){\line(0,-1){3}}
\put(3.1,3.1){\line(1,2){1.65}}
\put(3.1,3.1){\line(2,1){3.3}}
\put(3.1,3.1){\line(1,1){3}}
\put(4,3.5){$\bullet$}
\put(3.5,4){$\bullet$}
\put(4.5,4.5){$\bullet$}
\put(4.25,3.25){(2,1)}
\put(2.5,4){(1,2)}
\put(4.7,4.5){(3,3)}
\put(1.5,3){$\bullet$}
\put(3,1.5){$\bullet$}
\put(3.25,1.5){$(0,-3)$}
\put(1,3.3){$(-3,0)$}
\put(3,3){$\bullet$}
\put(2.1,2.6){$(0,0)$}
\end{picture}\\
\vspace{3mm}
Figure 1
\end{center}

One can check by calculating the semistable loci that this fan itself is the GIT fan. By Proposition \ref{movable} one also knows that the cone generated by (2,1) and (1,2) is the movable cone.

We can also investigate the smoothness of $X$. It was already proven by Bellamy \cite{Bel} that $V/G$ admits a symplectic resolution. We now try to prove the same thing using the Cox ring. To do this, let's consider two $\mathbb{C}^*$-actions on $V^*$ defined by $(x,y,z,w)\mapsto(tx,ty,tz,tw)$ and $(x,y,z,w)\mapsto(x,y,tz,tw)$ for $t\in\mathbb{C}^*$. Since these actions are compatible with the $G$-action on $V$, they induce the actions of $\mathbb{C}^*$ on $V/G$. Kaledin showed that these actions on $V/G$ lift to $X$ and the common fixed point sets consists of finitely many points \cite{Ka2}. In our case, these actions also lift to $\mathfrak{X}:=\mathrm{Spec}\,\mathrm{Cox}(X)$, and the fixed point set on $\mathfrak{X}$ is a single $T(=(\mathbb{C}^*)^2)$-orbit $$F=\{\psi_1=\cdots=\psi_{14}=\psi_{16}=\cdots=\psi_{19}
=\psi_{21}=0,\psi_{15}\ne0,\psi_{20}\ne0\}\subset\mathfrak{X}$$
where $\psi_i$ is the $i$-th generator of the Cox ring in the previous section regarded as the coordinate of $\mathfrak{X}$. Therefore, if the semistable locus of $\mathfrak{X}$ has a singular point, it must be in $F$. Note that this subset is contained in the semistable locus with respect to both of the two chambers in the movable cone. As we already know the explicit generators of the Cox ring, we can obtain their relations by a computer calculation, see Appendix \ref{App2}. Then the Jacobian criterion shows that $\mathfrak{X}$ is nonsingular at any point in $F$. One can check that each point of the semistable locus has a nontrivial stabilizer subgroup $T'\subset T$ of order 3 and that the quotient torus $T/T'$ acts freely on it. Therefore, by Luna's \'{e}tale slice theorem, one can conclude that $X$ is also smooth.

\section{(Non)smoothness of the minimal models of some symplectically imprimitive quotient singularities \label{6}}

In this section we investigate the smoothness of the minimal models for several cases. Now we are particularly interested in the symplectic cases. By Verbitsky's result, we only have to check the groups which are generated by symplectic reflections. Such groups are classified by Cohen \cite{C}. In his original paper he considered quaternion reflection groups rather than symplectic reflection groups, but one sees that these two kinds of groups can be identified.

To explain the classification, we prepare some terminologies. Let $V$ be a finite dimensional symplectic $\mathbb{C}$-vector space and let $\omega$ be its symplectic form. Let $Sp(V,\omega)$ (or simply $Sp(V)$) be the group of linear automorphisms of $V$ which preserve $\omega$ and let $G$ be a finite subgroup of $Sp(V)$. We say that the subgroup (or the representation) $G\subset Sp(V)$ is {\em irreducible} if there are no nontrivial decomposition of $V$ into $G$-invariant symplectic vector subspaces. Since every representation $G\subset Sp(V)$ is decomposed into irreducible representations, we will only consider irreducible ones from now on.

An irreducible $G$ is called {\em improper} if there is a $G$-invariant Lagrangian subspace $L$ of $V$ with respect to $\omega$ and otherwise we call $G$ {\em proper}. If $G$ is improper with $L\subset V$, the symplectic reflection group $G$ can be regarded as a complex reflection group via the natural inclusion $GL(L)\subset Sp(V)
$ \cite{V}. Complex reflection groups are classified by Shepherd-Todd \cite{ST} into three infinite families and 34 exceptional groups $G_4,\dots,G_{37}$.

Proper groups are also divided into two classes. We say that $G$ is {\em symplectically imprimitive} if there is a nontrivial decomposition
$$V=V_1\oplus\cdots\oplus V_k$$
into symplectic subspaces such that for any $g\in G$ and $i\in\{1,\dots,k\}$ there is $j$ such that $g(V_i)\subset V_j$. Otherwise $G$ is called {\em symplectically primitive}.

Bellamy and Schedler studied when quotient singularities by symplectically imprimitive groups have projective symplectic resolutions \cite{BS2}. They showed there that if $\mathrm{dim}\;V>4$, then $V/G$ has a symplectic resolution if and only if $G$ is the wreath product of a finite subgroup of $SL(2,\mathbb{C})$ and a symmetric group. Four dimensional (irreducible) proper symplectically imprimitive representations are classified up to conjugacy by Cohen and listed in the table I in \cite{C}. However one should note that Cohen's list is incomplete since it includes some improper groups and mutually conjugate groups as we will see later. We call the groups in the Cohen's list type (A), type (B),$\dots$,type (V) as in \cite{BS2}. Bellamy and Schedler also determined which $V/G$ has a symplectic resolution except 6 cases: type (G), type (K), type (P), type (Q), type (U), and type (V). The aim of this section is to complete their work by studying these remaining cases.

Let $V=\mathbb{C}^4$ and $\omega=dx\wedge dy+dz\wedge dw$ where $x,y,z$ and $w$ is the standard coordinate on $\mathbb{C}^4$. Then any of the 6 groups is of the following form
$$G(K,\alpha)=\bigcup_{i=1,2}\bigcup_{x\in K}\begin{pmatrix}x&0\\0&\alpha(x)\end{pmatrix}
\begin{pmatrix}0&1\\1&0\end{pmatrix}^i$$
where $K$ is a finite subgroup of $SL(2,\mathbb{C})$ and $\alpha\in\mathrm{Aut}(K)$ is an involution.

The six cases are listed as follows.

\vspace{3mm}
\noindent {\bf type (G)}$_{l,r}\,(l,r\in\mathbb{N}\text{ such that }r\le l,\,r\text{ is odd, and }l=\mathrm{gcd}(l,\frac{r+1}{2})\mathrm{gcd}(l,\frac{r-1}{2}))$:\\
$K=\langle g_1=\begin{pmatrix}\zeta&0\\0&\zeta^{-1}\end{pmatrix},
g_2=\begin{pmatrix}0&\sqrt{-1}\\\sqrt{-1}&0\end{pmatrix}\rangle\;(\zeta=\mathrm{exp}\,\frac{2\pi\sqrt{-1}}{2l})$
is a binary dihedral group and $\alpha$ is defined by $\alpha(g_1)=g_1^r,\alpha(g_2)=-g_2$.

\vspace{3mm}
\noindent {\bf type (K)}: $K=\langle g_1=\begin{pmatrix}\sqrt{-1}&0\\0&-\sqrt{-1}\end{pmatrix},
g_2=\frac{1}{\sqrt{2}}\begin{pmatrix}\zeta^5&\zeta^5\\\zeta^7&\zeta^3\end{pmatrix}\rangle\;(\;\zeta=\mathrm{exp}\,\frac{2\pi\sqrt{-1}}{8})$
is a binary tetrahedral group and $\alpha$ is defined by $\alpha(g_1)=\begin{pmatrix}0&-1\\1&0\end{pmatrix}$ and $\alpha(g_2)=\frac{1}{\sqrt{2}}\begin{pmatrix}\zeta^3&\zeta\\\zeta^3&\zeta^5\end{pmatrix}$.

\vspace{3mm}
\noindent {\bf type (P)}: $K=\langle g_1=\begin{pmatrix}\zeta&0\\0&\zeta^{-1}\end{pmatrix},
g_2=\frac{1}{\sqrt{2}}\begin{pmatrix}\zeta^5&\zeta^5\\\zeta^7&\zeta^3\end{pmatrix}\rangle
\;(\zeta=\mathrm{exp}\,\frac{2\pi\sqrt{-1}}{8}$)
is a binary octahedral group and $\alpha$ is defined by $\alpha(g_1)=g_1^{-1}$ and $\alpha(g_2)
=\frac{1}{\sqrt{2}}\begin{pmatrix}\zeta^3&\zeta^7\\\zeta^5&\zeta^5\end{pmatrix}$.

\vspace{3mm}
\noindent {\bf type (Q)}: $K=\langle g_1=\begin{pmatrix}\zeta&0\\0&\zeta^{-1}\end{pmatrix},
g_2=\frac{1}{\sqrt{2}}\begin{pmatrix}\zeta^5&\zeta^5\\\zeta^7&\zeta^3\end{pmatrix}\rangle
\;(\zeta=\mathrm{exp}\,\frac{2\pi\sqrt{-1}}{8}$)
is a binary octahedral group and $\alpha$ is defined by $\alpha(g_1)=-g_1$ and $\alpha(g_2)=g_2$.

\vspace{3mm}
\noindent {\bf type (U)}: $K=\langle g_1=\frac{1}{2}\begin{pmatrix}
\phi+\sqrt{-1}\phi^{-1}&1\\-1&\phi-\sqrt{-1}\phi^{-1}\end{pmatrix},
g_2=\frac{1}{\sqrt{2}}\begin{pmatrix}\zeta^5&\zeta^5\\\zeta^7&\zeta^3\end{pmatrix}\rangle
\;(\phi=\frac{1+\sqrt{5}}{2},\zeta=\mathrm{exp}\,\frac{2\pi\sqrt{-1}}{8}$)
is a binary icosahedral group and $\alpha$ is defined by\\ $\alpha(g_1)
=g_1^{-1}$ and
$\alpha(g_2)
=\frac{1}{\sqrt{2}}\begin{pmatrix}
\zeta^3&\zeta^7\\\zeta^5&\zeta^5\end{pmatrix}$.

\vspace{3mm}
\noindent {\bf type (V)}: $K=\langle g_1=\frac{1}{2}\begin{pmatrix}
\phi+\sqrt{-1}\phi^{-1}&1\\-1&\phi-\sqrt{-1}\phi^{-1}\end{pmatrix},
g_2=\begin{pmatrix}\sqrt{-1}&0\\0&-\sqrt{-1}\end{pmatrix}\rangle
\;(\phi=\frac{1+\sqrt{5}}{2}$)
is a binary icosahedral group and $\alpha$ is defined by\\ $\alpha(g_1)=-g_1g_2g_1^3$ and
$\alpha(g_2)
=\begin{pmatrix}0&\sqrt{-1}\\\sqrt{-1}&0\end{pmatrix}$.

\vspace{5mm}
Our main result in this section is the following.

\begin{thm}\label{classify}
Let $G$ be one of the 6 types above. Then the quotient singularity $V/G$ has a projective symplectic resolution if and only if $G$ is of type (G)$_{1,1}$.
\end{thm} 
{\em Proof.} We prove the claim by case-by-case analysis. First consider type (P) and type (U). These groups are improper groups (with respect to $\omega$). Indeed, one can easily check that the Lagrangian subspace $L=\{x-w=y-z=0\}$ of $V$ is preserved by the actions of the two groups. The corresponding complex reflection groups to type (P) and type (U) are $G_{13}$ and $G_{22}$ in the Shepherd-Todd classification \cite{ST} respectively. By Bellamy's result \cite{Bel}, we know that $V/G$ for each of $G_{13}$ and $G_{22}$ dose not have projective symplectic resolutions.

Next we consider type (K) and type (V). To cope with these groups, we consider the Cox rings of minimal models. Since direct computer calculations of the Cox rings could not be done in a reasonable amount of time, we adopt another approach.

Let $G$ denote the group of type (K) and $G'$ denote the group of type (J). Then $G'=\langle G,g_2=\begin{pmatrix}I_2&0\\0&-I_2\end{pmatrix}\rangle$ and $G$ is a normal subgroup of $G'$ of index 2. Let $g_1$ be a representative of the unique junior conjugacy class in $G$. The commutator groups $[G,G]$ and $[G',G']$ are the same, and we let $H$ denote this subgroup. By the results of section \ref{4}, the Cox ring of a minimal model $X$ for type (K) and that of a minimal model $X'$ for type (J) are realized as subrings of $R_1:=\mathbb{C}[V]^{H}[t_1^{\pm1}]$ and $R_2:=\mathbb{C}[V]^{H}[t_1^{\pm1},t_2^{\pm1}]$ respectively. Let
$$\psi_1=\phi_1t_1^{\nu_1(\phi_1)}t_2^{\nu_2(\phi_1)},\dots,\psi_k=\phi_kt_1^{\nu_1(\phi_k)}t_2^{\nu_2(\phi_k)}, T_1=t_1^{-2},T_2=t_2^{-2}\in R_2$$
be the generators of $\mathrm{Cox}(X')$ (see Proposition \ref{prop:generator}) which are homogeneous with respect to $G'/H(\cong\mathbb{Z}/4\mathbb{Z})$-action where $\phi_1\dots,\phi_k\in\mathbb{C}[V]^{H}$. By Proposition \ref{prop:generator}, we see that $\psi_1|_{t_2=1},\dots,\psi_k|_{t_2=1}, T_1\in R_1$ are generators of $\mathrm{Cox}(X)$. Since $\psi_i$'s are homogeneous with respect to $\langle g_2\rangle$-action, the $G'/G$-action on $V/G$ lifts to $\mathrm{Spec}\,\mathrm{Cox}(X)$. This action descends to one on the GIT quotient $X$ since the semistable locus is defined by the homogeneous elements by Proposition \ref{prop:chamber}. Note that the fixed point set of this action on $X$ is the common zero locus of $\psi_i|_{t_2=1}$'s such that $\nu_2(\phi_i)$ is odd.

By Proposition \ref{prop:chamber}, the GIT chambers on $\mathrm{Cl}(X')_{\mathbb{R}}\cong\mathbb{R}^2$ are described as in the following figure.

\begin{center}
\setlength\unitlength{1truecm}
\begin{picture}(6,6)(0,0)
\put(3.1,3.1){\line(-1,0){3}}
\put(3.1,3.1){\line(0,1){3}}
\put(3.1,3.1){\line(1,0){3}}
\put(3.1,3.1){\line(0,-1){3}}
\put(3.1,3.1){\line(1,5){0.7}}
\put(3.1,3.1){\line(2,5){1.2}}
\put(3.1,3.1){\line(3,1){3.3}}
\put(3.1,3.1){\line(3,2){3}}
\put(5,3.25){$C_1$}
\put(4,2.5){$D$}
\put(4.5,1.5){$C_2$}
\put(4.5,5.6){$\cdot$}
\put(4.85,5.4){$\cdot$}
\put(5.1,5.2){$\cdot$}
\put(5.3,4.9){$\cdot$}
\put(4,3){$\bullet$}
\put(1.5,3){$\bullet$}
\put(3,1.5){$\bullet$}
\put(3.25,1.5){$E_2$}
\put(1,3.3){$E_1$}
\put(3,3){$\bullet$}
\put(2.5,2.6){$\mathbf{O}$}
\end{picture}\\
\vspace{3mm}
Figure 2
\end{center}

The divisor $D=-\frac{1}{2}E_1$ in the figure is obtained, for example, as the degree of the element in the Cox ring associated to the semi-invariant $x^{12}-33x^8y^4-33x^4y^8+y^{12}\in\mathbb{C}[V]^{H}$. By Proposition \ref{movable}, one sees that the half line $\mathbb{R}_{\ge0}D$ lies on the boundary of the movable cone.
Note that the GIT quotient of the spectrum $\mathfrak{X}':=\mathrm{Spec}\,\mathrm{Cox}(X')$ with respect to $D\in\mathrm{Cl}(X')_{\mathbb{R}}$ is the quotient of $X$ by the $\langle g_2\rangle(\cong\mathbb{Z}/2\mathbb{Z})$-action. Note also that the GIT quotient of $\mathfrak{X}'$ with respect to the open chamber $C_2$ in the figure is the same as $\mathfrak{X}'/\!/_D T$ (cf. \cite[1.11]{HK}). We may assume that $X'$ is the minimal model which corresponds to the open chamber $C_1$ in the figure. The semistable loci on $\mathfrak{X}'$ with respect to $C_1,C_2$ and $D$ have following inclusions:
$$\mathfrak{X}'^{C_1-ss}\subset\mathfrak{X}'^{D-ss}\supset\mathfrak{X}'^{C_2-ss}.$$
By the definition of a stable point, we also see that $\mathfrak{X}'^{D-s}=\mathfrak{X}'^{C_1-s}\cap\mathfrak{X}'^{C_2-s}$. Recall that the morphism $\pi:\mathfrak{X}'/\!/_{C_1} T\to\mathfrak{X}'/\!/_D T=X/\langle g_2\rangle$ and the isomorphism $\mathfrak{X}'/\!/_{C_2} T\to\mathfrak{X}'/\!/_D T$ of GIT quotients are induced from the inclusions of the semistable loci on $\mathfrak{X}'$. Therefore we see that $\pi$ is an isomorphism on the image of $\mathfrak{X}'^{D-s}$ in $\mathfrak{X}'/\!/_{C_1} T$
. One can directly check by Proposition \ref{prop:chamber} that $\mathfrak{X}'^{C_1-s}\backslash \mathfrak{X}'^{D-s}=\mathfrak{X}'^{C_1-s}\cap \{T_2=0\}$ and $\mathfrak{X}'^{C_2-s}\backslash \mathfrak{X}'^{D-s}=\mathfrak{X}'^{C_2-s}\cap \{\psi_i=0|\nu_2(\phi_i)>0\}$. Therefore $\pi$ contracts the unique irreducible exceptional divisor $E_2$ defined by $\{T_2=0\}$ onto the set $F\subset X/\langle g_2\rangle$ which is defined by $\psi_i$'s such that $\nu_2(\phi_i)>0$.

Now we assume that $X$ is smooth. Since the $\langle g_2\rangle$-action is symplectic, the singularities of $X/\langle g_2\rangle$ is analytically locally isomorphic to $\mathbb{C}^2\times(\mathbb{C}^2/\{\pm1\})$ or $\mathbb{C}^4/\{\pm1\}$. Since the isolated singularity $\mathbb{C}^4/\{\pm1\}$ is already terminal, the singularity of $X/\langle g\rangle$ along $F$ is isomorphic to $\mathbb{C}^2\times(\mathbb{C}^2/\{\pm1\})$. (Note that $X/\langle g_2\rangle$ may have singularities of type $\mathbb{C}^4/\{\pm1\}$ outside $F$.) Thus the blowing-up $X'$ of $X/\langle g_2\rangle$ along $F$ must be smooth in a neighborhood of $E_2$. Therefore, in order to prove that $X$ is singular, it suffices to show that $X'$ has singularities in $E_2$.

In \cite{BS2} the authors consider the minimal resolution $Y$ of $\mathbb{C}^2/\{\pm1\}$ and show that some minimal model $X''\to V/G'$ factors through $(Y\times Y)/H'$ for $H'=G'/\langle g_2,\begin{pmatrix}-I_2&0\\0&I_2\end{pmatrix}\rangle$. Note that the exceptional divisor of $(Y\times Y)/H'\to V/G'$ is associated to the symplectic reflection $g_2$. One can directly check that $(Y\times Y)/H'$ has a singularity which is isolated (and hence terminal) using the argument in \cite[5.3]{BS2}. Thus $X''$ has a singular point in the exceptional divisor associated to $g_2$. Since $X'$ and $X''$ are connected by a sequence of Mukai flops \cite[Thm. 1.2]{WW}, $X'$ also has a singular point in $E_2$. Therefore $X$ is singular. Note that the minimal model $X$ is unique since $\mathrm{Cl}(X)_{\mathbb{R}}$ is 1-dimensional.

For type (V), we can use the exactly same argument as one for type (K) by replacing type (J) by type (T).

For type (Q), the group $G$ is $Sp(V,\omega)$-conjugate to another group in Cohen's list. Indeed, one can easily check that the matrix $$g=\frac{\sqrt{-1}}{2}\begin{pmatrix}\zeta&-\zeta^3&-\zeta&\zeta^3\\
\zeta&\zeta^3&-\zeta&-\zeta^3\\
1&-1&1&-1\\
-1&-1&-1&-1\end{pmatrix}$$
is in $Sp(V,\omega)$ and that the group of type (J) is the $g$-conjugate of $G$ where $\zeta=\mathrm{exp}\,\frac{2\pi\sqrt{-1}}{8}$. Therefore $V/G$ does not admit projective symplectic resolutions by \cite{BS2}.

Finally we treat type (G)$_{l,r}$. When $r$ is 1, one can check that $G$ preserves the Lagrangian subspace $L=\{x+\sqrt{-1}w=y-\sqrt{-1}z=0\}$ of $V$ and thus $G$ is improper (with respect to $\omega$). The corresponding complex reflection is $G(2l,l,2)$ in the Shepherd-Todd classification \cite{ST}. By the result of Bellamy \cite{Bel}, we know that $V/G$ with $G=G(2l,l,2)$ admits a projective symplectic resolution if and only if $l=1$.
 
When $r\ne1$, we use the similar method as type (K) and type (V). The group $G'=\langle G,g=\begin{pmatrix}I_2&0\\0&-I_2\end{pmatrix}\rangle$ coincides with the group of type (B) (resp. type (F)) in Cohen's list if $l$ is even (resp. odd). By the same argument above using the description of Cox rings, one sees that $g$-action on $V/G$ lifts to its minimal model $X$ and that a minimal model $X'$ of $V/G'$ is obtained as the blowing-up of $X/\langle g\rangle$. Similarly to the cases type (K) and type (V), it suffices to show that $X'$ has a singular point in the exceptional divisor $E$ which corresponds to the symplectic reflection $g$.

Consider the same $Y$ as above and $H'=G'/\langle g,\begin{pmatrix}-I_2&0\\0&I_2\end{pmatrix}\rangle$. In this case $Y\times Y$ has no isolated fixed points by $H'$-action but one can find a point $x\in Y\times Y$ such that $\mathrm{Stab}_{H'}(x)$ is not generated by symplectic reflections by using the argument in \cite[5.3]{BS2}. Let $\mathrm{Stab}_{H'}(x)^\circ$ be the normal subgroup of $\mathrm{Stab}_{H'}(x)$ generated by symplectic reflections. One can also find $x$ such that there is $h\in\mathrm{Stab}_{H'}(x)\backslash\mathrm{Stab}_{H'}(x)^\circ$ whose action preserves the exceptional divisors of $Y\times Y\to V/\langle g,\begin{pmatrix}-I_2&0\\0&I_2\end{pmatrix}\rangle$. This implies that the minimal model $X'$ has a singular point in $E$.
\qed

\section{Appendix \label{7}}
\subsection{}\label{App1}
In this subsection we give a concrete method to perform the algorithm. We will usually need computer calculations in practice.

Let $G\subset SL(V)$ be a finite subgroup and let $g_1,\dots,g_m$ be a complete system of representatives of the conjugacy classes of the junior elements in $G$. Assume that we are given the generators $\phi_1,\dots,\phi_k$ of the invariant ring $\mathbb{C}[V]^{[G,G]}$ which are homogeneous with respect to $Ab(G)$-action.

Let $I$ be the kernel of $\alpha:\mathbb{C}[X_1,\dots,X_k]\to\mathbb{C}[V]^{[G,G]}$ (see section \ref{4}). Recall that for each $i\in\{1,\dots,m\}$, the variable $X_j$ has degree $\nu_i(\phi_j)$. Let $\mathbb{C}[X_1,\dots,X_k,t]$ be the new graded ring  with one more variable $t$ whose degree is $-1$. The ideal $\mathrm{min}_i(I)\;(i=1,\dots,m)$ is calculated by taking the following steps:\\
1. Consider the ideal generated by the image of $I$ by the inclusion $\mathbb{C}[X_1,\dots,X_k]\hookrightarrow$\\
\hspace{5mm}$\mathbb{C}[X_1,\dots,X_k,t]$, and let this ideal also denote $I$ by abuse of notation.\\
2. Homogenize the generators of $I$ with respect to the variable $t$, and let $I_i$ be the ideal\\
\hspace{5mm}generated by these elements.\\
3. Compute the saturation $\tilde{I}_i:=\bigcup_{l=0}^\infty I_i:(t^l)$ of $I_i$ with respect to $t$.\\
4. Evaluate $t=0$ in $\tilde{I}_i$.\\
Then the resulting ideal (regarded as the ideal of $\mathbb{C}[X_1,\dots,X_k]$) is $\mathrm{min}_i(I)$. Note that just homogenizing the generators is not enough and the saturation is necessary in general.

In order to obtain $I_A\subset R_A$ for a subset $A\subset\{1,\dots,m\}$, one should consider $R_A=\mathbb{C}[X_1,\dots,X_k,\{t_i\}_{i\in A}]$ and perform the steps from 1. to 3. for each $t_i\,(i\in A)$. 

Next let's consider $\mathrm{min}_i(J)\;(i=1,\dots,m)$. Let $g\in G$ and let $r$ be the order of $g$ in $Ab(G)$. As each $\phi_j$ is homogeneous, there is an integer $0\le a_j<r$ such that $g$ acts on $\phi_j$ by multiplication of $\mathrm{exp}\,\frac{2\pi\sqrt{-1}a_j}{r}$. Let $\mathbb{C}[X_1,\dots,X_k,s]$ be the graded polynomial ring where $\deg(X_j)=a_i$ and $\deg(s)=1$. Let $J$ be the kernel of $\beta_i$ (see  section \ref{4}). We can calculate the ideal $J_g$ generated by homogeneous elements of $J$ with respect to $g$ by taking following steps:\\
1.  Consider the ideal generated by the image of $J$ by the inclusion $\mathbb{C}[X_1,\dots,X_k]\hookrightarrow$\\
\hspace{5mm}$\mathbb{C}[X_1,\dots,X_k,s]$, and let this ideal also denote $J$ by abuse of notation.\\
2. Homogenize the generators of $J$ with respect to the variable $s$ and let $J_i$ be the ideal\\
\hspace{5mm}generated by these elements.\\
3. Compute the saturation $\tilde{J}_i:=\bigcup_{l=0}^\infty J_i:(s^l)$ of $J_i$ with respect to $s$.\\
Then the preimage of $\tilde{J}_i+(s^r-1)$ by the inclusion $\mathbb{C}[X_1,\dots,X_k]\hookrightarrow\mathbb{C}[X_1,\dots,X_k,s]$ is $J_g$. By repeating the same procedures over $g$'s which generate $Ab(G)$, we finally obtain $\mathrm{min}_i(J)$.

\subsection{}\label{App2}
In this subsection we give a relatively easy way of calculation of the relations of the generators of the Cox ring.

By the algorithm, we know that the generators of the Cox ring of the minimal model $X$ of $V/G$ is of the following form:
$$\psi_1:=\phi_1\prod_{i=1}^m t_1^{\nu_i(\phi_1)},\dots,\psi_k:=\phi_1\prod_{i=1}^m t_k^{\nu_i(\phi_k)},T_1=t_1^{-r_1},\dots,T_m:=t_m^{-r_m}$$
where $\phi_i$'s are the homogeneous generators of $\mathbb{C}[V]^{[G,G]}$, $g_1,\dots,g_m$ are the representatives of the junior elements in $G$ and $r_i:=\sharp\langle g_i\rangle$. Assume that we are already given the ideal $I\subset\mathbb{C}[X_1,\dots,X_k]$ of the relations of $\phi_i$'s. Then the ideal $\tilde{I}\subset\mathbb{C}[X_1,\dots,X_k,Y_1,\dots,Y_m]$ of the relations of $\psi_i$'s and $T_i$'s are calculated as follows:\\
1. Compute $I_{\{1,\dots,m\}}\subset R_{\{1,\dots,m\}}=\mathbb{C}[X_1,\dots,X_k,t_1,\dots,t_m]$ (see \ref{App1}).\\
2. Replace every $t_i^{r_i}$ by $Y_i$ in the $Ab(G)^\vee$-homogeneous generators of  $I_{\{1,\dots,m\}}$ for each $i$. (This is possible since homogeneity implies that $t_i$'s appear only with powers of multiples of $r_i$.)\\
The resulting ideal is $\tilde{I}$.


\vspace{0.2cm}

\begin{center}
Department of Mathematics, Faculty of Science, Kyoto University, Japan 

ryo-yama@math.kyoto-u.ac.jp
\end{center}

\end{document}